\pdfoutput=1

\documentclass[12pt]{amsart}
\usepackage[margin=1in]{geometry} 

\usepackage{amsmath,amssymb,amsthm,bm,colonequals,graphicx,mathrsfs,mathtools,microtype,stmaryrd}
\usepackage[shortlabels]{enumitem}
\usepackage[colorinlistoftodos]{todonotes}

\theoremstyle{plain}

\newtheorem{theorem}{Theorem}[section] 

\newtheorem{lemma}[theorem]{Lemma} 

\newtheorem{proposition}[theorem]{Proposition} 

\newtheorem{proposition-definition}[theorem]{Proposition-Definition} 

\newtheorem{corollary}[theorem]{Corollary} 

\theoremstyle{definition}

\newtheorem{definition}[theorem]{Definition} 

\newtheorem{problem}[theorem]{Problem} 

\newtheorem{question}[theorem]{Question} 

\theoremstyle{remark}

\newtheorem{remark}[theorem]{Remark} 

\newcommand{\Aff}{\mathbb{A}} 

\newcommand{\CC}{\mathbb{C}}
\newcommand{\EE}{\mathbb{E}}
\newcommand{\FF}{\mathbb{F}}

\newcommand{\PP}{\mathbb{P}}
\newcommand{\QQ}{\mathbb{Q}}

\newcommand{\ZZ}{\mathbb{Z}}

\newcommand{\rank}{\operatorname{rank}}
\newcommand{\im}{\operatorname{im}}

\newcommand{\Span}{\operatorname{Span}}

\newcommand{\abs}[1]{\lvert #1 \rvert}
\newcommand{\card}[1]{\lvert #1 \rvert}
\newcommand{\norm}[1]{\lVert #1 \rVert}

\newcommand{\eps}{\epsilon}

\newcommand{\Spec}{\operatorname{Spec}}

\newcommand{\inject}{\operatorname{\hookrightarrow}}

\newcommand{\disc}{\operatorname{disc}}

\newcommand{\map}{\operatorname}

\newcommand{\mcal}{\mathcal}
\newcommand{\mf}{\mathfrak}

\newcommand{\ol}{\overline}

\newcommand{\defeq}{\colonequals}
\newcommand{\eqdef}{\equalscolon}
\newcommand{\maps}{\colon}

\newcommand{\belongs}{\subseteq}
\newcommand{\contains}{\supseteq}

\newcommand{\set}[1]{\{#1\}}

\newcommand{\grad}{\nabla}

\usepackage[pagebackref=true]{hyperref}
\usepackage[alphabetic,initials]{amsrefs}

\begin{document}

\title{Dichotomous point counts over finite fields}
\author{Victor Y. Wang}
\address{Department of Mathematics, Princeton University, Princeton, NJ, USA}
\address{Courant Institute of Mathematical Sciences, New York University, New York, NY, USA}
\email{vywang@alum.mit.edu}
\date{}

\begin{abstract}
We establish a near dichotomy between randomness and structure
for the point counts of arbitrary projective cubic threefolds over finite fields.
Certain ``special'' subvarieties, not unlike those in the Manin conjectures, dominate.
We also prove new general results for projective hypersurfaces.
Our work continues a line of inquiry initiated by Hooley.
\end{abstract}

\maketitle

\setcounter{tocdepth}{3}

\section{Introduction}


Let $X$ be a projective scheme over a finite field $\FF_q$, and consider the point counts $\#X(\FF_{q^r})$ for integers $r\geq 1$.
Say that $X$ \emph{satisfies} (resp.~\emph{fails}) ``square-root cancellation'' if
\begin{equation}
\label{INEQ:|E|-good}
\sup_{r\geq 1}
\left(\frac{\abs{\#X(\FF_{q^r}) - \#\PP^{\dim X}_{\FF_q}(\FF_{q^r})}}{q^{r(\dim X)/2}}\right)
< \infty
\end{equation}
holds (resp.~fails).
(For an ``effective'' equivalent of \eqref{INEQ:|E|-good},
see e.g.~Proposition~\ref{PROP:amplification-for-projective-varieties}.)
The square-root normalization in \eqref{INEQ:|E|-good} is most natural for special classes of $X$ such as hypersurfaces and complete intersections;
for other $X$, there may be other, more natural normalizations.

Now suppose $X$ is a projective complete intersection over $\FF_q$.
By \cite{hooley1991number}*{Theorem~2} (a result building on Deligne's resolution of the Weil conjectures),
\begin{equation}
\label{INEQ:Hooley1991-implication}
\sup_{r\geq 1}
\left(\frac{\abs{\#X(\FF_{q^r}) - \#\PP^{\dim X}_{\FF_q}(\FF_{q^r})}}{q^{r(1+\dim(X_{\map{sing}})+\dim{X})/2}}\right)
< \infty.
\end{equation}
Here $\dim(X_{\map{sing}})$ denotes the dimension of the singular locus of $X$,
and the empty set has dimension $-1$.
If $X$ is smooth, then \eqref{INEQ:|E|-good} and \eqref{INEQ:Hooley1991-implication} coincide.
For singular $X$,
there are examples (e.g.~iterated cones over a smooth hypersurface of degree $\ge 3$) where \eqref{INEQ:Hooley1991-implication} is essentially optimal (in the sense that the exponent $1+\dim(X_{\map{sing}})+\dim{X}$ in the denominator cannot be reduced),
but \cite{hooley1991number}*{second paragraph after Theorem~2} makes the following points (among others):
\begin{itemize}
    \item \eqref{INEQ:|E|-good} fails for \emph{some} singular choices of $X$.
    
    \item It would be nice to have a ``satisfactory criterion'' to determine when \eqref{INEQ:|E|-good} fails.
    
    \item Some ``scanty evidence'' suggests that \eqref{INEQ:|E|-good} rarely fails to the extent that \eqref{INEQ:Hooley1991-implication} allows.
\end{itemize}
Yet \eqref{INEQ:Hooley1991-implication} itself is silent on these three points.
Our Theorems~\ref{THM:main-general-result-informal}, \ref{THM:main-special-result-informal}, \ref{THM:main-diagonal-theorem}, and~\ref{THM:concrete-moment-and-super-level-bounds}(2)--(3) below go beyond \eqref{INEQ:Hooley1991-implication} in this respect,
for certain $X=V_{\bm{c}}$ (isomorphic to hypersurfaces).

Throughout this paper, we let $V_U(f)_{/R}$ denote the closed subscheme of $U_R$ cut out by $f=0$;
and we let $V_U(f_1,\dots,f_r)_{/R}$ denote the (scheme-theoretic) intersection $\bigcap_{i=1}^{r} V_U(f_i)_{/R}$.
If the base ring $R$
or the ambient object $U$
is clear from context, we may omit it.

Given a base field $k$,
an integer $m\geq 3$,
and a nonconstant homogeneous polynomial $F\in k[\bm{x}]=k[x_1,\dots,x_m]$,
let $V$ denote the $(m-2)$-dimensional projective hypersurface
\begin{equation}
\label{EQN:projective-V(F)}
V_{\PP^{m-1}}(F)_{/k}
= \set{[\bm{x}]\in \PP^{m-1}_k: F(\bm{x}) = 0},
\end{equation}
and for $\bm{c}=(c_1,\dots,c_m)\in k^m\setminus \set{\bm{0}}$ let $V_{\bm{c}}$ (or $V_{[\bm{c}]}$) denote the hyperplane section
\begin{equation}
\label{EQN:projective-V(F,c.x)}
V_{\PP^{m-1}}(F,\bm{c}\cdot\bm{x})_{/k}
= \set{[\bm{x}]\in \PP^{m-1}_k: F(\bm{x}) = \bm{c}\cdot\bm{x} = 0}
\end{equation}
(which is of dimension $m-3$
if, for instance, $V$ is smooth and $\deg{F}\geq 2$).
We will repeatedly use this setup with $k,m,F,V,V_{\bm{c}}$; call it the \emph{Main Setup} for convenience.

The set notation on the right-hand side of \eqref{EQN:projective-V(F)} and \eqref{EQN:projective-V(F,c.x)} is meant to be suggestive but slightly informal (as $V, V_{\bm{c}}$ are really schemes, not sets).
We abbreviate projective coordinates using vector notation;
thus, for instance, if $\bm{x}=(x_1,\dots,x_m)\in \ol{k}^m\setminus \set{0}$, then $[\bm{x}]$ denotes the point $[x_1:\dots:x_m]\in \PP^{m-1}_k(\ol{k})$.
Also, $\bm{c}\cdot\bm{x} \defeq c_1x_1+\dots+c_mx_m$ denotes the usual dot product.

\begin{theorem}
\label{THM:main-general-result-informal}
Adopt the Main Setup.
Suppose $k$ is finite,
$V$ is smooth,
$d\defeq \deg{F}\geq 3$,
and $\map{char}(k)\nmid d(d-1)$.
Then there is a closed subscheme $B\belongs \Aff^m_k$ of codimension $\geq 2$, with $\bm{0}\in B$, such that $V_{\bm{c}}$ satisfies ``square-root cancellation'' for all $\bm{c}\in k^m\setminus B(k)$.
One can take $B\defeq V_{\Aff^m}(\disc(F,\bm{c}), \partial\disc(F,\bm{c})/\partial c_1, \dots, \partial\disc(F,\bm{c})/\partial c_m)_{/k}$, where $\disc$ is defined as in \S\ref{SEC:background-on-discriminants}.
\end{theorem}


Theorem~\ref{THM:main-general-result-informal} does not extend to $\deg{F}=2$ in general (see the beginning of \S\ref{SEC:establishing-near-dichotomy-for-quadrics-and-low-dimensional-cubics}),
but its truth for $\deg{F}=3$ is significant for dependent work of the author on (diagonal) cubic forms (see \cite{wang2022thesis}*{Chapters~7--8}).
We will prove Theorem~\ref{THM:main-general-result-informal} in \S\ref{SEC:general-codimension-two-work}, mainly by combining results on discriminants, duality, and $\ell$-adic cohomology (especially \cite{lindner2020hypersurfaces}*{Theorem~1.2}, but also work of Deligne, Katz, Skorobogatov, et al.~to be discussed and repackaged in \S\ref{SEC:general-background}).
On the other hand, Sawin has informed us that one could perhaps prove a less explicit version of Theorem~\ref{THM:main-general-result-informal} using the perversity-based strategy of \cite{grimmelt2021representation}*{proof of Lemma~3.1}.

Theorem~\ref{THM:main-general-result-informal} is fairly general.
But if one believes that ``randomness'' should increase with $d,m$, then yet stronger statements should hold for sufficiently large $d,m$.
In this paper, we formulate and prove \emph{optimal} statements in special cases.
The quadratic case is more or less understood (in odd characteristic, at the very least), so we focus on cubics.
Theorem~\ref{THM:main-special-result-informal} suggests, and Theorem~\ref{THM:main-diagonal-theorem} proves, that Theorem~\ref{THM:main-general-result-informal} is far from the full truth in general.

\begin{theorem}
\label{THM:main-special-result-informal}
In the Main Setup,
assume $k$ is finite,
$V$ is smooth,
$\deg{F}=3$,
and $m\in \set{4,6}$.
Let $\bm{c}\in k^m\setminus \set{\bm{0}}$.
Then $V_{\bm{c}}$ fails ``square-root cancellation'' \emph{if and only if} the base change $V_{\bm{c}}\times_k \ol{k}$ contains
either an $\frac{m-2}{2}$-plane, or a singular cubic $2$-scroll (as defined in \S\ref{SEC:establishing-near-dichotomy-for-quadrics-and-low-dimensional-cubics}), in $\PP^{m-1}_{\ol{k}}$.
\end{theorem}

\begin{theorem}
\label{THM:main-diagonal-theorem}
In the Main Setup,
assume $k$ is finite,
$m\in \set{4,6}$,
and $F = F_1x_1^3+\dots+F_mx_m^3$ for some $F_1,\dots,F_m\in \ZZ\setminus \set{0}$.
Let $\bm{c}\in k^m\setminus \set{\bm{0}}$.
If $\map{char}(k)$ is sufficiently large in terms of $F_1,\dots,F_m$,
then $V_{\bm{c}}$ fails ``square-root cancellation'' \emph{if and only if} ``$c_i^3/F_i=c_j^3/F_j$ in pairs''
(i.e.~there exist $m/2$ pairwise disjoint sets $\set{i,j}\belongs \ZZ$ with $1\leq i<j\leq m$ and $c_i^3/F_i = c_j^3/F_j$).
\end{theorem}

Theorems~\ref{THM:main-special-result-informal}--\ref{THM:main-diagonal-theorem}
(proven in \S\S\ref{SEC:establishing-near-dichotomy-for-quadrics-and-low-dimensional-cubics}--\ref{SEC:Fermat-cubic-fourfold-analysis} using classical geometry, on top of some of the modern geometry behind Theorem~\ref{THM:main-general-result-informal})
are new for $m=6$,
and are closely connected to work of the author on special subvarieties in Manin-type conjectures (see \cite{wang2022thesis}*{Chapter~6}).
This connection suggests the following vague question;
see \cite{wang2022thesis}*{Question~5.1.1} for more details.

\begin{question}
\label{QUES:correlate-special-subvarieties-Manin-vs-F_q}
To what extent do special subvarieties in the Manin conjectures correlate with special subvarieties in the sense suggested by Theorem~\ref{THM:main-special-result-informal}?
\end{question}

One cannot expect special subvarieties to completely explain all failures of ``square-root cancellation''; for example, a general cone $X$ of sufficiently high degree may well provide a counterexample.
However, one might still hope that special subvarieties could explain ``most'' failures in a suitable algebraic or statistical sense.
See Remark~\ref{RMK:analytic-auxiliary-polynomial-complexity-direction-of-inquiry} for an alternative hope.

There is certainly room for interesting theoretical and numerical testing and experimentation; we have only scratched the surface of possibilities.
Full numerical testing, however, may be hard.
One can usually test a statement of the form ($\bm{c}$ lies in a given algebraic set)$\Rightarrow$($V_{\bm{c}}$ fails ``square-root-cancellation''), but less often the converse.
The data in \cite{github-singular-cubic-threefold-2021}*{{\tt Singular...~tests 4.0}--{\tt 4.2.txt}} gives moderate, but not strong, numerical evidence for Theorem~\ref{THM:main-diagonal-theorem}.
Even with thousands of tests (each of which may take millions or billions of operations),
with testing restricted to the locus $\disc(F,\bm{c})=0$ (say),
one can only hope to reliably detect phenomena over $\FF_{101}$ (say) of codimension $2$ or $3$ in $\Aff^m$.
Based on data alone, it would be too easy to conjecture a false version of the general Theorem~\ref{THM:main-special-result-informal}.
(Singular cubic $2$-scrolls can be ``rarer'' than planes; see \S\ref{SEC:Fermat-cubic-fourfold-analysis}.)

We now measure Theorems~\ref{THM:main-general-result-informal}, \ref{THM:main-special-result-informal}, and~\ref{THM:main-diagonal-theorem} against a concrete problem:


\begin{problem}
\label{PROB:moments-of-failures-of-sqrt-cancellation}
Fix $G\in \ZZ[\bm{x}]=\ZZ[x_1,\dots,x_m]$, homogeneous of degree $d\geq 2$ in $m\geq 3$ variables,
with nonzero discriminant
(so that $V_{\PP^{m-1}}(G)_{/\CC}$ is smooth of dimension $m-2$).
Given a prime $p$ and integer $r\geq 1$, let $q\defeq p^r$.
Given $\bm{c}=(c_1,\dots,c_m)\in \FF_q^m$, let
\begin{equation}
	E_{\bm{c}}(q)\defeq \#V_{\PP^{m-1}}(G, \bm{c}\cdot\bm{x})_{/\FF_q}(\FF_q) - \#\PP^{m-3}_{\FF_q}(\FF_q).
\end{equation}
For real $\sigma\geq 0$, consider the moments
\begin{equation}
	\label{DEFN:E_c-moment}
	M_G(q,\sigma)\defeq \sum_{\bm{c}\in \FF_q^m\setminus \set{\bm{0}}} \abs{E_{\bm{c}}(q) / q^{(m-3)/2}}^{\sigma}
\end{equation}
and the associated ``exponents''
\begin{equation}
	\label{DEFN:E_c-moment-exponent}
	e_G(\sigma)\defeq \limsup_{p\to \infty}
	\left(\sup_{r\geq 1} \left(\log_q(1 + M_G(q,\sigma))\right)\right).
\end{equation}
As $G,d,m,\sigma$ vary, estimate $e_G(\sigma)$.
\end{problem}

The moments $M_G(q,\sigma)$ are closely related to the super-level sets
\begin{equation}
	\label{DEFN:lambda-bad-count}
	S_G(q,\eps)\defeq \set{\bm{c}\in \FF_q^m\setminus \set{\bm{0}}: \abs{E_{\bm{c}}(q)} \geq q^{\eps + (m-3)/2}}
\end{equation}
(for real $\eps\geq 0$, say).
Both $M_G(q,\sigma)$ and $\#S_G(q,\eps)$ concern the large values of $\abs{E_{\bm{c}}(q)}$.

The first part of the following theorem is essentially classical (due to Deligne, Zak, Hooley, et al.), but the rest is new (to the author's knowledge).

\begin{theorem}
\label{THM:concrete-moment-and-super-level-bounds}
Let $G,d,m,p,r,q$ be as in Problem~\ref{PROB:moments-of-failures-of-sqrt-cancellation}.
Fix a real $\eps>0$.
\begin{enumerate}
    \item There exists a real constant $C=C(G,d,m,\eps)>0$ such that if $p\geq C$, then $\#S_G(q,\eps)\leq C q^{m-1}$ and $M_G(q,2)\leq C q^m$.
    In particular, $e_G(2)\leq m$.
    
    \item Suppose $d\geq 3$.
    There exists $C=C(G,d,m,\eps)>0$ such that if $p\geq C$, then $\#S_G(q,\eps)\leq C q^{m-2}$ and $M_G(q,4)\leq C q^m$.
    In particular, $e_G(4)\leq m$.
    
    \item Suppose $(m,d) =  (6,3)$ and $G$ is diagonal, i.e.~of the form $G_1x_1^d+\dots+G_mx_m^d$.
    There exists $C=C(G,d,m,\eps)>0$ such that if $p\geq C$, then $\#S_G(q,\eps)\leq C q^{m-3}$ and $M_G(q,6)\leq C q^m$.
    In particular, $e_G(6)\leq m$.
\end{enumerate}
\end{theorem}

We will prove Theorem~\ref{THM:concrete-moment-and-super-level-bounds} in Appendix~\ref{SEC:moment-calculations}, using Theorems~\ref{THM:main-general-result-informal} and~\ref{THM:main-diagonal-theorem}.

The problem of bounding $\#S_G(q,\eps)$ is close in spirit to \cite{lindner2020hypersurfaces}*{Theorem~1.3 and the line after},
but Lindner studies a different aspect (with $q,m$ fixed and $d\to \infty$), works with a different family (a universal family of hypersurfaces), and works not with $E_{\bm{c}}(q)$ but with a related cohomological quantity.
One could likely extend most of our main results to universal families of hypersurfaces; or reformulate some of them in terms of cohomological ``defect'' (in the sense of \cite{lindner2020hypersurfaces}), via Proposition~\ref{PROP:amplification-for-projective-varieties} below; or both.
But we focus on the $V_{\bm{c}}$'s and their point counts, since they play---in a fashion building on \cite{hooley1986HasseWeil}---a central role in \cite{wang2022thesis}.

It seems likely that the diagonality assumption on $G$ in Theorem~\ref{THM:concrete-moment-and-super-level-bounds}(3) is unnecessary.
We do not know how to prove this, but see \S\ref{SEC:Fermat-cubic-fourfold-analysis} for further discussion.

Also, could one evaluate $M_G(q,\sigma)$ asymptotically, or  relate $M_G(q,\sigma)$ to other moments in number theory (e.g.~moments of $L$-functions)?
For $M_G(q,2)$, see e.g.~Appendix~\ref{SEC:moment-calculations}, whose main goal is to prove the following corollary to Theorem~\ref{THM:main-general-result-informal}.
When $2\mid m$, Corollary~\ref{COR:smooth-locus-moment-calculations} strengthens (in two of the three error terms) a consequence of the Deligne--Katz equidistribution theorem.
(The case $2\nmid m$ may or may not involve ``exceptional monodromy'', which could complicate an attempt to use Deligne--Katz.)
Thus Corollary~\ref{COR:smooth-locus-moment-calculations} suggests that there may be a deeper connection between monodromy and Problem~\ref{PROB:moments-of-failures-of-sqrt-cancellation} than one might at first expect.

\begin{corollary}
[Proven in Appendix~\ref{SEC:moment-calculations}]
\label{COR:smooth-locus-moment-calculations}
Let $G,d,m,E_{\bm{c}}(q)$ be as in Problem~\ref{PROB:moments-of-failures-of-sqrt-cancellation}.
If $d,m\geq 3$, then the following hold uniformly over $p$ (with implied constants depending only on $G$):
\begin{enumerate}
    \item $\EE_{\bm{c}\in \FF_p^m}[E_{\bm{c}}(p)\cdot \bm{1}_{p\nmid \disc(G,\bm{c})}]
    = p^{-1}\cdot E(V_{\PP^{m-1}}(G)_{/\FF_p}) + p^{-1}\cdot O(p^{(m-3)/2})$.
    
    \item $\EE_{\bm{c}\in \FF_p^m}[E_{\bm{c}}(p^2)\cdot \bm{1}_{p\nmid \disc(G,\bm{c})}]
    = (1+O(p^{-1/2}\bm{1}_{m=3})+O(p^{-1}))\cdot p^{m-3}$.
    
    \item $\EE_{\bm{c}\in \FF_p^m}[E_{\bm{c}}(p)^2\cdot \bm{1}_{p\nmid \disc(G,\bm{c})}]
    = (1+O(p^{-1/2}\bm{1}_{m=3})+O(p^{-1}))\cdot p^{m-3}$.
\end{enumerate}
Here $\EE_{\bm{c}\in S}[f(\bm{c})] \defeq \card{S}^{-1} \sum_{\bm{c}\in S} f(\bm{c})$.
(For the notation $\bm{1}_A, E(X), O(\cdot)$, see Definitions~\ref{DEFN:analytic-and-combinatorial-notation} and~\ref{DEFN:algebro-geometric-notation}.
For the notation $\disc(G,\bm{c})$, see \S\ref{SEC:background-on-discriminants}.)
\end{corollary}

The rest of the paper is organized as follows:
\S\ref{SEC:general-background} sets up some general conventions and background,
\S\ref{SEC:background-on-discriminants} gives background on discriminants,
\S\ref{SEC:general-codimension-two-work} proves Theorem~\ref{THM:main-general-result-informal} (and a bit more),
\S\ref{SEC:establishing-near-dichotomy-for-quadrics-and-low-dimensional-cubics} proves Theorem~\ref{THM:main-special-result-informal} (and a bit more),
\S\ref{SEC:Fermat-cubic-fourfold-analysis} proves Theorem~\ref{THM:main-diagonal-theorem} (and discusses the diagonality assumption in Theorem~\ref{THM:concrete-moment-and-super-level-bounds}(3)),
and Appendix~\ref{SEC:moment-calculations} proves Theorem~\ref{THM:concrete-moment-and-super-level-bounds} and Corollary~\ref{COR:smooth-locus-moment-calculations}.

\section{General background}
\label{SEC:general-background}

The main goals of this section are to provide Definitions~\ref{DEFN:analytic-and-combinatorial-notation}, \ref{DEFN:algebro-geometric-notation}, and~\ref{DEFN:error-goodness-for-projective-X/F_q} below (all used throughout the paper), and to state and prove Proposition~\ref{PROP:amplification-for-projective-varieties} below.

\begin{definition}
[Combinatorial and analytic conventions]
\label{DEFN:analytic-and-combinatorial-notation}
We let $[r]\defeq \set{1,2,\dots,r}$ for positive integers $r$.
We let $\bm{1}_A$ denote the \emph{indicator value} of an event $A$;
i.e.~we let $\bm{1}_A\defeq 1$ if $A$ holds, and $\bm{1}_A\defeq 0$ otherwise.
We write $f\ll g$ to mean $\abs{f}\leq Cg$ for some real constant $C>0$.
We let $O(g)$ denote a quantity that is $\ll g$.
\end{definition}

In this paper, we work in terms of schemes for convenience in proofs and applications.
We only use the term ``variety'' informally (mainly to emphasize connections with literature), and in these instances the concrete classical definition ``integral, locally closed subscheme of a projective space over a field'' would be perfectly adequate.

\begin{definition}
[Geometric conventions]
\label{DEFN:algebro-geometric-notation}
Let $k$ be a base field.
A \emph{scheme over $k$} (or \emph{$k$-scheme} for short) is a scheme over $\Spec{k}$;
a \emph{curve} (resp.~\emph{surface}) is a $k$-scheme of dimension $1$ (resp.~$2$).
An \emph{embedded projective $k$-scheme} is a projective $k$-scheme realized as a closed subscheme of some $\PP^n_k$.
A \emph{hypersurface} in $\PP^n_k$ is an embedded projective $k$-scheme equal to $V_{\PP^n}(P)_{/k}$ for some nonconstant homogeneous $P\in k[x_1,\dots,x_{n+1}]$.
An intersection of $r$ hypersurfaces in $\PP^n_k$ is called a \emph{complete intersection} if it has dimension $n-r$.

In the context of projective $k$-schemes,
let $\dim(\emptyset)\defeq -1$ and $\PP^{-1}_k\defeq \emptyset$.
If $X$ is a $k$-scheme and $L$ is a field extension of $k$, let $X(L)$ denote the set of $L$-points on $X$, as in any standard reference (e.g.~\cite{poonen2017rational}*{\S2.3.2}).
If $X$ is a projective $k$-scheme and $k$ is finite, let
\begin{equation}
\label{DEFN:E(X)}
E(X)\defeq \#X(k)-\#\PP^{\dim{X}}_k(k).
\end{equation}
(The base field $k$ of $X$ is essential to this definition.
In particular, if $k'/k$ is a finite extension, then $E(X_{k'}) = \#X(k')-\#\PP^{\dim{X}}_k(k')$.
Here $X_{k'}$ denotes the \emph{base change} $X\times_k k'$ of $X$ to $k'$.)

Say a scheme is of \emph{pure dimension $d$} if all of its irreducible components are of dimension $d$.
For a finite-type $k$-scheme $X$ of pure dimension $d\geq 0$,
let $X_{\map{sing}}$ denote the \emph{singular subscheme} of $X$,
i.e.~the closed subscheme of $X$ defined by the $d$th Fitting ideal of the cotangent sheaf $\Omega_{X/k}$ (or informally, ``by the Jacobian criterion''), following \cite{stacks-project}*{\href{https://stacks.math.columbia.edu/tag/0C3H}{Tag~0C3H}}.
For a scheme $Y$, we let $\abs{Y}$ denote the underlying topological space;
thus, for instance, $\abs{X_{\map{sing}}}$ denotes a topological space,
while $\card{X_{\map{sing}}(\ol{k})} = \#X_{\map{sing}}(\ol{k})$ denotes an integer.
\end{definition}

We now recall a useful principle of Zak (often used to verify hypotheses on $X_{\map{sing}}$):

\begin{theorem}
[Zak; see e.g.~\cite{hooley1991number}*{Katz's Appendix, Theorem~2}]
\label{THM:Zak's-principle}
Let $k,m,F,V,V_{\bm{c}}$ be as in the Main Setup.
If $V$ is smooth, then $\dim((V_{\bm{c}})_{\map{sing}})\leq 0$ for all $\bm{c}\in k^m\setminus \set{\bm{0}}$.
\end{theorem}

It is worth noting that in most of the paper, the scheme structure on $X_{\map{sing}}$ is not essential.
A notable exception is \S\ref{SEC:Fermat-cubic-fourfold-analysis}, where we use the following general technical fact.
(It would be nice to know how much the hypotheses can be weakened; cf.~\cite{MO211153intersect_components}.)

\begin{proposition}
\label{PROP:lower-bound-on-singular-scheme-of-a-reducible-variety}
Let $X_1,X_2$ be locally closed subschemes of $\PP^n$ of pure dimension $d$ over an algebraically closed field $K$, where $d,n\geq 1$.
Let $X\defeq X_1\cup X_2$.
Assume the following:
\begin{enumerate}
    \item $\dim(X_1\cap X_2) = 0$;
    and
    \item for each $x\in X_1\cap X_2$, there exist a locally closed subscheme $Y$ of $\PP^n$ of pure dimension $2d$, and an open neighborhood $U$ of $x$ in $\PP^n$,
    such that $X_1\cap U$ and $X_2\cap U$ are Cohen--Macaulay, $Y\cap U$ is smooth, and $Y\cap U\contains X\cap U$.
\end{enumerate}
Then $X_{\map{sing}}\contains X_1\cap X_2$.
\end{proposition}

\begin{proof}
The statement is local, so we may assume $X_1\cap X_2$ is supported on a singleton $\set{x}$.
Now let $Y,U$ be as in hypothesis~(2); by shrinking $U$ if necessary, we may assume that $U$ is affine and that $X_1,X_2,Y$ are closed subschemes of $U$.
Say $U = \Spec{R}$, and let $I_1,I_2,J\belongs R$ be the ideals defining $X_1,X_2,Y$, respectively.
Then $X_1,X_2$ are Cohen--Macaulay, $Y$ is regular, and $J\belongs I_1\cap I_2$.
By \cite{SpeyerMO49261product_ideal} and \cite{serre2000local}*{Proposition~11 in \S{IV.B.1}, and Corollary to Theorem~4 in \S{V.B.6}},
it follows that $(I_1/J)\cap (I_2/J) = (I_1/J)\cdot (I_2/J)$ in $R/J$ (cf.~\cite{DaoMO49299product_ideal}).
So if $f\in I_1\cap I_2$, then $f\equiv h\bmod{I_1I_2}$ for some $h\in J$, whence $Df\equiv Dh\bmod{(I_1,I_2)}$ for all derivations $D\maps R\to R$.
But $\Omega_{Y/K}$ is locally free of rank $2d\geq d+1$, so the $d$th Fitting ideal of $\Omega_{Y/K}$ is $0$.
It follows that any $(n-d)\times (n-d)$ matrix of the form $(D_j f_i)_{i,j}$ (where $f_i\in I_1\cap I_2$) has determinant $\equiv 0\bmod{(I_1,I_2)}$.
Thus $V_U(I_1\cap I_2)_{\map{sing}}\contains V_U(I_1,I_2)$, i.e.~$X_{\map{sing}}\contains X_1\cap X_2$.
\end{proof}

The following definition places ineq.~\eqref{INEQ:|E|-good} (from the introduction) in a broader network of ``error-goodness'' concepts (concepts related to one another by Proposition~\ref{PROP:amplification-for-projective-varieties}).

\begin{definition}
\label{DEFN:error-goodness-for-projective-X/F_q}
Let $k$ denote an arbitrary finite field, and $X$ a projective $k$-scheme.
\begin{enumerate}
    \item Given $f\in \set{\abs{E}, +E, -E}$ (where $E$ is the function $E(-)$ defined by eq.~\eqref{DEFN:E(X)}),
    say $X$ is \emph{$f$-good} if there exists a real $C>0$ such that for all finite extensions $k'/k$, we have
    \begin{equation*}
    f(X_{k'})\leq C\card{k'}^{(\dim{X})/2}.
    \end{equation*}
    If $C$ is admissible, say $X$ is \emph{$f$-good with constant $C$}.
    Say $X$ is \emph{$f$-bad} if it is not $f$-good.
    
    \item Given a property \emph{blah} in (1), say $X$ is \emph{potentially} (resp.~\emph{stably}) \emph{blah} if $X_{k'}$ is blah for some (resp.~for every) finite extension $k'/k$.
\end{enumerate}
\end{definition}

In the setting of Definition~\ref{DEFN:error-goodness-for-projective-X/F_q},
it is routine to prove the following lemma (for use in \S\ref{SEC:establishing-near-dichotomy-for-quadrics-and-low-dimensional-cubics}):

\begin{lemma}
\label{LEM:potential-and-stable-logic}
Let $K/k$ be a finite extension.
Then $X$ is potentially $f$-good if and only if $X_K$ is potentially $f$-good.
Similarly, $X$ is stably $f$-bad if and only if $X_K$ is stably $f$-bad.







\end{lemma}

For the rest of \S\ref{SEC:general-background}, let $k$ denote a finite field.
We now state a useful amplification-type result.
We will prove it using general foundational results due to Deligne, Katz, and others.

\begin{proposition}
\label{PROP:amplification-for-projective-varieties}
Let $X$ be the intersection of $r$ hypersurfaces in $\PP^n_k$, each of degree $\leq d$, where $n,r,d\geq 1$.
Consider the following six conditions:
\begin{enumerate}
    \item\label{ITEM:effectively-error-good} $X$ is $\abs{E}$-good with constant $18(3+rd)^{n+1}2^r$;
    \item\label{ITEM:error-good} $X$ is $\abs{E}$-good;
    \item\label{ITEM:potentially-error-good} $X$ is potentially $\abs{E}$-good;
    \item\label{ITEM:potentially-signed-error-good} $X$ is potentially $(-1)^{1+\dim{X}}E$-good;
    \item\label{ITEM:no-defect-in-1+dim} in the notation of Definition~\ref{DEFN:nontrivial-or-error-relevant-Frob-eigenvalues-for-projective-X/F_q} below,
    $\dim H^{1+\dim{X}}(X) = \dim H^{1+\dim{X}}(\PP^{\dim{X}}_k)$ holds for every prime $\ell\nmid \card{k}$;
    and
    \item\label{ITEM:absolutely-error-good} $X$ is absolutely $\abs{E}$-good,
    in the sense of Definition~\ref{DEFN:absolute-error-goodness-for-projective-X/F_q} below.
\end{enumerate}
In general,
(\ref{ITEM:effectively-error-good})--(\ref{ITEM:error-good}) are equivalent,
(\ref{ITEM:error-good}) implies (\ref{ITEM:potentially-error-good}),
(\ref{ITEM:potentially-error-good}) implies (\ref{ITEM:potentially-signed-error-good}),
and (\ref{ITEM:absolutely-error-good}) implies (\ref{ITEM:effectively-error-good})--(\ref{ITEM:error-good}).
If $\dim{X} = n-r$
and $\dim(X_{\map{sing}})\leq 0$,
then (\ref{ITEM:effectively-error-good})--(\ref{ITEM:absolutely-error-good}) are equivalent.
\end{proposition}

\begin{remark}
In particular, a projective complete intersection $Y/k$ with $\dim(Y_{\map{sing}})\leq 0$ is $\abs{E}$-good if and only if it is potentially $\abs{E}$-good.
It would be nice to have this more generally.
\end{remark}


We now build up to a proof of Proposition~\ref{PROP:amplification-for-projective-varieties}.
In view of the Grothendieck--Lefschetz trace formula for $E(X)$, it is natural to first make some cohomological definitions.

\begin{definition}
\label{DEFN:nontrivial-or-error-relevant-Frob-eigenvalues-for-projective-X/F_q}
Fix a prime $\ell\nmid \card{k}$.
Let $Y,X$ denote arbitrary projective $k$-schemes.
\begin{enumerate}
    \item For integers $i\geq 0$, define $H^i(Y)\defeq H^i(Y\times_k \ol{k}, \QQ_\ell)$ using $\ell$-adic cohomology with $\QQ_\ell$-coefficients.
    
    \item For each $i\geq 0$, let $\mcal{E}^i(Y)$ denote the multiset of (geometric) Frobenius eigenvalues on $H^i(Y)$.
    
    \item Let $\mcal{E}^i_{\triangle}(X,\PP)\defeq (\mcal{E}^i(X)\cup \mcal{E}^i(\PP^{\dim{X}}_k))\setminus (\mcal{E}^i(X)\cap \mcal{E}^i(\PP^{\dim{X}}_k))$ for $i\geq 0$.
    In other words, if $\alpha\in \ol{\QQ}_\ell$ has multiplicities $j_1,j_2\geq 0$ in $\mcal{E}^i(X),\mcal{E}^i(\PP^{\dim{X}}_k)$, let it have multiplicity $\abs{j_1-j_2}$ in $\mcal{E}^i_{\triangle}(X,\PP)$.
    
    \item Define $\mcal{E}_{\triangle}(X,\PP)\defeq \sum_{i\geq 0}\mcal{E}^i_{\triangle}(X,\PP)$ by ``summing multiplicities'' over $i$.
\end{enumerate}
\end{definition}

\begin{remark}
Deligne's purity theorem implies that each $\mcal{E}^i(Y)$ above consists of $\card{k}$-Weil numbers $\alpha\in \ol{\QQ}_\ell$ of weight $w(\alpha)\leq i$.
(See e.g.~\cite{kiehl2001weil}*{Remark~I.7.2 and Theorem~I.9.3(2)} for a precise textbook reference in English.)
\end{remark}

\begin{definition}
\label{DEFN:absolute-error-goodness-for-projective-X/F_q}
Let $X$ be a projective $k$-scheme.
Say $X$ is \emph{absolutely $\abs{E}$-good} if for every $\ell\nmid \card{k}$, all $\alpha\in \mcal{E}_{\triangle}(X,\PP)$ have weight $\leq \dim{X}$.
\end{definition}

\begin{remark}
In principle, eigenvalues of weight $\geq 1+\dim{X}$ could cancel out in the trace formula for $E(X)$.
But ``absolute $\abs{E}$-goodness'' forbids such subtleties.
\end{remark}

The following statement combines \cite{hooley1991number}*{Katz's Appendix, assertion~(2) in the proof of Theorem~1} and \cite{ghorpade2008etale}*{first sentence of Remark 3.5}.
Both are important for us, but the latter seems to appear without proof in \cite{ghorpade2008etale}, so we sketch some.

\begin{theorem}
[Katz, Skorobogatov, et al.]
\label{THM:general-perversity-result}
Fix $\ell\nmid \card{k}$.
Fix integers $n,N\geq 0$ with $n\geq N+1$,
and a complete intersection $X\belongs \PP^n_k$ with $\dim{X} = N$.
Let $D\defeq \dim(X_{\map{sing}})$.
If $i\in \ZZ$, then
(1) if $i\geq N+D+2$, then $\mcal{E}^i_{\triangle}(X,\PP) = \emptyset$;
and (2) if $i = N+D+1$, then $\mcal{E}^i(\PP^{N}_k)\belongs \mcal{E}^i(X)$.
\end{theorem}

\begin{remark}
By \cite{ghorpade2008etale}*{Proposition~3.2}, one could add the following to Theorem~\ref{THM:general-perversity-result}:
(3) if $i = N$, then $\mcal{E}^i(\PP^{N}_k)\belongs \mcal{E}^i(X)$;
and (4) if $0\leq i\leq N-1$, then $\mcal{E}^i_{\triangle}(X,\PP) = \emptyset$.
\end{remark}

\begin{proof}
[Proof sketch for Theorem~\ref{THM:general-perversity-result}]
When $X$ is a hypersurface section of a smooth projective complete intersection $Y/k$ with $\dim{Y}\geq 2$, Theorem~\ref{THM:general-perversity-result} follows from \cite{skorobogatov1992exponential}*{Corollary~2.2, up to Veronese embedding}, Theorem~\ref{THM:Zak's-principle}, ``Betti comparison'' for $Y$, and the geometric irreducibility of $Y$.
In general, one can prove Theorem~\ref{THM:general-perversity-result} either inductively or directly.
Suppose $D\leq N-2$ (and $N\geq 1$); the result follows from \cite{poonen2017rational}*{Corollary~7.5.21} otherwise.

\emph{Inductive proof.}
Induct on $n-N$ using \cite{skorobogatov1992exponential}*{Corollary~2.2}, \cite{ghorpade2008etale}*{Lemma~1.1(ii)}, and \cite{ghorpade2008etale}*{proofs of Theorem~2.4 and Proposition~2.5, up to Veronese embeddings}.

\emph{Direct proof.}
Claim~(1) follows from Katz \cite{hooley1991number}*{loc. cit.}.
For (2), follow Katz, but use the full strength of \cite{grothendieck1972groupes}*{Deligne's Expos\'{e}~I, Corollaire~4.3}.
\end{proof}

The following standard result is also essential in the proof of Proposition~\ref{PROP:amplification-for-projective-varieties}:

\begin{lemma}
[Real amplification]
\label{LEM:standard-real-amplification-lemma}
Fix an integer $l\geq 0$.
Let $\beta_1,\dots,\beta_l \in \set{z\in \CC: \abs{z}\geq 1}$.
Then $\limsup_{n\to \infty}\Re(\beta_1^{dn}+\dots+\beta_l^{dn})\geq l$ holds for every integer $d\geq 1$.
\end{lemma}

\begin{proof}
We may assume $l\geq 1$ and $d=1$.
Then, Dirichlet's approximation theorem implies that $\limsup_{n\to \infty}\Re({\cdots})$ is $+\infty$ if $\max(\abs{\beta_1},\dots,\abs{\beta_l}) > 1$, and $l$ otherwise.
\end{proof}

The present paper overall could do without the following purity result, which we nonetheless include to keep Proposition~\ref{PROP:amplification-for-projective-varieties} as clean as possible.

\begin{proposition}
[Bertini, Ghorpade--Lachaud, et al.]
\label{PROP:purity-in-degrees-N+D+1-and-higher}
In the setting of Theorem~\ref{THM:general-perversity-result}, if $i\geq N+D+1$, then all $\alpha\in \mcal{E}^i(X)$ have weight equal to $i$.
\end{proposition}

\begin{proof}
[Proof sketch]
Suppose $D\leq N-2$; the result follows from \cite{poonen2017rational}*{Corollary~7.5.21} otherwise.
Then $N\geq 1$, and $X_{\ol{k}}$ is regular in codimension one, so $X_{\ol{k}}$ is integral (by e.g.~\cite{ghorpade2008etale}*{\S3, par.~2}).
By the Bertini-based \cite{ghorpade2008etale}*{proof of Proposition~3.3(i)--(ii)},
one can find a finite extension $k'/k$, and a smooth projective complete intersection $Y\belongs \PP^n_{k'}$ of dimension $N-D-1$, admitting a surjective $G_{k'}$-equivariant map
\begin{equation*}
H^{i-2D-2}(Y_{\ol{k}}, \QQ_\ell(-D-1))
\to H^i(X_{\ol{k}}, \QQ_\ell)
\end{equation*}
for each $i\geq N+D+1$.
The left-hand side is pure of weight $i$ as a $G_{k'}$-module (by Deligne),
so the right-hand side is too
(as a $G_{k'}$-module, and thus also as a $G_k$-module).
\end{proof}

\begin{proof}
[Proof of Proposition~\ref{PROP:amplification-for-projective-varieties}]
The case $X=\emptyset$ is trivial, so assume $X\neq \emptyset$.
Let $N\defeq \dim{X}\geq 0$.

(\ref{ITEM:effectively-error-good})$\Rightarrow$(\ref{ITEM:error-good}), (\ref{ITEM:error-good})$\Rightarrow$(\ref{ITEM:potentially-error-good}), and (\ref{ITEM:potentially-error-good})$\Rightarrow$(\ref{ITEM:potentially-signed-error-good}):
Obvious.

(\ref{ITEM:error-good})$\Rightarrow$(\ref{ITEM:effectively-error-good}):
For each $m\geq 1$, let $k_m$ denote a degree-$m$ extension of $k$.
Assume (\ref{ITEM:error-good});
then the formal power series $\sum_{m\geq 1} E(X_{k_m})T^m\in \ZZ[[T]]$ defines a holomorphic function on the disk $\abs{T} < \card{k}^{-N/2}$.
But by the Grothendieck--Lefschetz trace formula,
$\sum_{m\geq 1} E(X_{k_m})T^m$ simplifies to
$\sum_{\beta\in S} M(\beta)\cdot \beta T(1-\beta T)^{-1}\in \CC(T)$ (a rational expression)
for some finite set $S\subset \CC$
and some function $M\maps S\to \ZZ\setminus \set{0}$ with
\begin{equation*}
\norm{M}_1
\defeq \sum_{\beta\in S} \abs{M(\beta)}
\leq \sum_{i}\sum_{Y\in \set{X, \PP^N}}\dim H^i(Y).
\end{equation*}
Thus in particular,
$\abs{\beta}\leq \card{k}^{N/2}$ for all $\beta\in S$,
whence $\abs{E(X_{k_m})}\leq \norm{M}_1\cdot (\card{k}^{N/2})^m$ for all $m\geq 1$.
But if $Y\in \set{X, \PP^N}$, then $\sum_{i}\dim H^i(Y)\leq 9(3+rd)^{n+1}2^r$, by the second inequality of \cite{katz2001sums}*{Corollary of Theorem~3}.
So (\ref{ITEM:effectively-error-good}) holds.

(\ref{ITEM:absolutely-error-good})$\Rightarrow$(\ref{ITEM:error-good}):
Assume (\ref{ITEM:absolutely-error-good}).
Fix $\ell\nmid \card{k}$.
Then the trace formula, Definition~\ref{DEFN:absolute-error-goodness-for-projective-X/F_q}, and finiteness of $\dim H^i(X), \dim H^i(\PP^N)$ give (\ref{ITEM:error-good}).

For the remainder of the proof, assume $\dim{X} = n-r$ and $\dim(X_{\map{sing}})\leq 0$.

(\ref{ITEM:potentially-signed-error-good})$\Rightarrow$(\ref{ITEM:absolutely-error-good}):
In general, $\dim(X_{\map{sing}})\leq 0$ and Theorem~\ref{THM:general-perversity-result} imply that the multiset $\set{\alpha\in \mcal{E}_{\triangle}(X,\PP): \map{weight}(\alpha)\geq 1+N}$ is a sub-multiset of $\mcal{E}^{1+N}(X)\setminus \mcal{E}^{1+N}(\PP^N_k)$.
So by the trace formula and Lemma~\ref{LEM:standard-real-amplification-lemma}, (\ref{ITEM:potentially-signed-error-good}) implies (\ref{ITEM:absolutely-error-good}).

(\ref{ITEM:no-defect-in-1+dim})$\Leftrightarrow$(\ref{ITEM:absolutely-error-good}):
If (\ref{ITEM:no-defect-in-1+dim}) holds,
then by Theorem~\ref{THM:general-perversity-result}, $\sum_{i\geq 1+N}\mcal{E}^i_{\triangle}(X,\PP) = \emptyset$ for all $\ell\nmid \card{k}$;
so (\ref{ITEM:absolutely-error-good}) holds.
If (\ref{ITEM:absolutely-error-good}) holds,
then by Proposition~\ref{PROP:purity-in-degrees-N+D+1-and-higher}, $\mcal{E}^{1+N}_{\triangle}(X,\PP) = \emptyset$ for all $\ell\nmid \card{k}$;
so (\ref{ITEM:no-defect-in-1+dim}) holds.
\end{proof}

\section{Background on discriminants}
\label{SEC:background-on-discriminants}

Let $d\geq 2$ and $m\geq 3$ be integers.
For a homogeneous polynomial $P$ of degree $d$ in $x_1,\dots,x_m$, let $P_{\bm{\alpha}}$ denote the coefficient of $\bm{x}^{\bm{\alpha}}=x_1^{\alpha_1}\cdots x_m^{\alpha_m}$ in $P$.
The \emph{discriminant} $\disc(P,\bm{c})$ is a polynomial in $\ZZ[\set{P_{\bm{\alpha}}},\bm{c}]=\ZZ[\set{P_{\bm{\alpha}}},c_1,\dots,c_m]$ characterized by certain properties.
Let us briefly recall the details.
First recall that a polynomial $Q$ over a field is said to be \emph{square-free} if $Q$ is not divisible by the square of any nonconstant polynomial.

\begin{proposition-definition}
[Cf.~\cite{terakado2018determinant}*{\S1.1}]
\label{PROPDEFN:discriminant-of-hyperplane-sections}
There exists a primitive homogeneous polynomial $D\in \ZZ[\set{P_{\bm{\alpha}}},\bm{c}]\setminus \set{0}$,
square-free over $\QQ$,
such that given $P,\bm{c}/\overline{\QQ}$, we have $D(P,\bm{c})\neq 0$ if and only if $V_{\PP^{m-1}}(P,\bm{c}\cdot\bm{x})_{/\ol{\QQ}}$ is smooth of dimension $m-3$.
Let $\disc\defeq D$.
Then $\disc$ is unique up to sign,
irreducible over $\ol{\QQ}$,
and bi-homogeneous of bi-degree $((m-1)(d-1)^{m-2}, d(d-1)^{m-2})$
in $(\set{P_{\bm{\alpha}}},\bm{c})$.
Also, over an arbitrary base ring $R$,
a given intersection $V_{\PP^{m-1}}(P,\bm{c}\cdot\bm{x})_{/R}$ is smooth of dimension $m-3$ if and only if $\disc(P,\bm{c})\in R^\times$.
\end{proposition-definition}

\begin{proof}
The existence, uniqueness, and irreducibility follow immediately from \cite{terakado2018determinant}*{paragraph before Proposition~1.2} (which quotes several results from \cite{benoist2012degres}).
For the degree claim, see \cite{benoist2012degres}*{Th\'{e}or\`{e}me~1.3}.
For the final statement, see \cite{terakado2018determinant}*{Proposition~1.3}.
\end{proof}

By e.g.~\cite{saito2012discriminant}*{\S2.1}, one can analogously define $\disc(G)$ for homogeneous $G$ of degree $d\geq 2$ in $m-1\geq 2$ variables; $\disc(-)$ is homogeneous of degree $(m-1)(d-1)^{m-2}$.
Although $\disc(G)$ is perhaps more familiar than $\disc(P,\bm{c})$, there is a close connection between the two:

\begin{proposition}
\label{PROP:d,1-vs-d-discriminant}
In $\ZZ[\set{P_{\bm{\alpha}}},c_1/c_m,\dots,c_{m-1}/c_m,c_m]$,
the elements
\begin{equation*}
\disc(P,\bm{c})
\quad\textnormal{and}\quad
\disc(P\vert_{x_m=-(c_1x_1+\dots+c_{m-1}x_{m-1})/c_m})\cdot c_m^{d(d-1)^{m-2}}
\end{equation*}
coincide up to sign.
In particular, both actually lie in $\ZZ[\set{P_{\bm{\alpha}}},\bm{c}]$.
\end{proposition}

\begin{proof}
Call the elements $L,R$, respectively.
Let $\Aff\defeq\Spec(k[\set{P_{\bm{\alpha}}},\bm{c}])$ over $k\defeq\overline{\QQ}$, and let $U\defeq \Aff \setminus V_{\Aff}(c_m)$.
Then $\abs{V_U(L)} = \abs{V_U(R)}$, by definition of $L,R$.
So if we choose $d\geq 0$ with $c_m^dR\in k[\Aff]$, then $\abs{V_{\Aff}(c_m^dR)}\belongs \abs{V_{\Aff}(c_m\cdot L)}$.
Since $c_m,L$ are \emph{prime} in $k[\Aff]$, it follows (by Hilbert's nullstellensatz) that $c_m^dR\in k^\times\cdot c_m^e\cdot L^f$ for some integers $e,f\geq 0$.

But in the variables $P_{\bm{\alpha}}$, each of $L,R$ is homogeneous of degree $(m-1)(d-1)^{m-2}$.
So $f=1$.
Yet in $c_1,\dots,c_m$, each of $L,R$ is homogeneous of degree $d(d-1)^{m-2}$.
Thus $e=d$, and $R\in k^\times\cdot L$.
In particular, $R\in\ZZ[\set{P_{\bm{\alpha}}},\bm{c}/c_m,c_m]\cap k[\set{P_{\bm{\alpha}}},\bm{c}]=\ZZ[\set{P_{\bm{\alpha}}},\bm{c}]$.
Here $R$ is primitive, since for the ``universal'' homogeneous polynomial $G=G(x_1,\dots,x_{m-1})$ of degree $d$ in $x_1,\dots,x_{m-1}$, we have $R\vert_{(P,c_m)=(G,1)} = \disc(G)\in \ZZ[\set{G_{\bm{\alpha}}}]$, which is primitive by definition.
But $L\in \ZZ[\set{P_{\bm{\alpha}}},\bm{c}]$ is itself primitive by definition.
So $R = \pm L$.
\end{proof}

If one restricts attention to \emph{diagonal} $P$, say $P_1x_1^3+\dots+P_mx_m^3$ (with $d=3$),
there is a nice explicit formula for $\disc(P,\bm{c})$.
(This generalizes to other degrees, but to avoid a mess of notation we focus on the cubic case.)
We record the (classical) details here for the interested reader.
Recall that the Fermat cubic $x_1^3+\dots+x_{m-1}^3$ has discriminant $\pm 3^{\mf{e}_m}$,
where
\begin{equation*}
\mf{e}_m
\defeq \frac{(-1)^{m-1} - 2^{m-1}}{3}
+ (m-1)2^{m-2}
= 3\cdot\bm{1}_{m=3}
+ 9\cdot\bm{1}_{m=4}
+ 27\cdot\bm{1}_{m=5}
+ 69\cdot\bm{1}_{m=6}
+\cdots
\end{equation*}
(for $m=3$, see \cite{saito2012discriminant}*{\S5.1};
for $m=4$, see \cite{saito2012discriminant}*{\S5.3};
in general, see \cite{gelfand2008discriminants}*{p.~434, (1.8) and Proposition~1.7}).

\begin{definition}
\label{DEFN:correct-disc-Delta(F,c)-definition-for-diagonal-F}
Given a \emph{diagonal} cubic form $P=P_1x_1^3+\dots+P_mx_m^3$,
let
\begin{equation*}
\Delta^\textnormal{diag}(P,\bm{c})
\defeq 3^{\mf{e}_m}
(P_1\cdots P_m)^{2^{m-2}}
\prod_{\bm{\eps}}
\left[\eps_1(c_1^3/P_1)^{1/2}
+ \eps_2(c_2^3/P_2)^{1/2}
+ \dots
+ \eps_m(c_m^3/P_m)^{1/2}\right],
\end{equation*}
where $\bm{\eps}=(\eps_1,\dots,\eps_m)$ ranges over $\set{1}\times\set{\pm1}^{m-1}$.
Here the right side represents an element of $\ZZ[\set{P_i},\bm{c}]$,
factored in $\ZZ[\set{P_i},\bm{c}][\set{(c_i^3/P_i)^{1/2}}]$.
\end{definition}

\begin{proposition}
\label{PROP:diag-disc-factorization}
If $P$ denotes the ``universal'' \emph{diagonal} cubic form $P_1x_1^3+\dots+P_mx_m^3$,
then the equality
$\disc(P,\bm{c})
= \pm \Delta^\textnormal{diag}(P,\bm{c})$
holds in $\ZZ[\set{P_i},\bm{c}]$.
\end{proposition}

\begin{proof}
At the level of $\ol{\QQ}$-points,
the polynomials $\disc(P,\bm{c}),\Delta^\textnormal{diag}(P,\bm{c})\in\QQ[\set{P_i},\bm{c}]$
have the same zero locus,
$\set{(P,\bm{c}): V(P,\bm{c}\cdot\bm{x})\;\textnormal{is not smooth of dimension $m-3$}}$.
(For proof, use the Jacobian criterion to compute the singularities of the projective system of equations $P=\bm{c}\cdot\bm{x}=0$.
The case $P_1\cdots P_m=0$ requires some care;
one must prove $\Delta^\textnormal{diag}(P,\bm{c})=0$ when
$P_i=c_i=0$ for some $i\in[m]$,
and also when $P_i=P_j=0$ for some distinct $i,j\in [m]$.)

Also, each of $\disc(P,\bm{c}), \Delta^\textnormal{diag}(P,\bm{c})$ has total degree $(m+2)2^{m-2}$,
and (by Lemma~\ref{LEM:bi-homogeneous-factored-polynomial-is-irreducible} below) $\Delta^\textnormal{diag}(P,\bm{c})$ is prime in $\QQ[\set{P_i},\bm{c}]$.
Hence $\disc(P,\bm{c}), \Delta^\textnormal{diag}(P,\bm{c})$ coincide up to a factor in $\QQ^\times$.
To finish,
note that equality (up to sign) holds
when we specialize $(P,\bm{c})$ to $(x_1^3+\dots+x_{m-1}^3, (0,\dots,0,1))$:
by Proposition~\ref{PROP:d,1-vs-d-discriminant},
$\disc(P,\bm{c})$ becomes
$\pm \disc(P\vert_{x_m=0})
= \pm 3^{\mf{e}_m}$,
while by direct computation,
$\Delta^\textnormal{diag}(P,\bm{c})$ becomes
``$\pm 3^{\mf{e}_m}P_m^{2^{m-2}}(c_m^3/P_m)^{2^{m-2}}$'' (i.e.~$\pm 3^{\mf{e}_m}$).
\end{proof}


\begin{lemma}
\label{LEM:bi-homogeneous-factored-polynomial-is-irreducible}
$\Delta^\textnormal{diag}(P,\bm{c})$ is irreducible in $\ol{\QQ}[\set{P_i},\bm{c}]$ (assuming $m\geq 3$ as above).
\end{lemma}

\begin{proof}
Let $k\defeq\ol{\QQ}$
and $A\defeq k[\set{P_i},\set{c_i}]$.
Consider the multi-quadratic ring extension $B\defeq A[\set{E_i},\set{d_i}]/(E_i^2-P_i,d_i^2-c_i)\cong_k k[\set{E_i},\set{d_i}]$.
Here $A,B$ are unique factorization domains,
and $B$ is ``Galois'' over $A$,
so it suffices to show that the element
\begin{equation*}
P
\defeq E_1\cdots E_m
\left[\eps_1d_1^3/E_1
+ \dots
+ \eps_md_m^3/E_m\right]
= \eps_1d_1^3E_2\cdots E_m
+ \dots
+ \eps_md_m^3E_1\cdots E_{m-1}
\in B
\end{equation*}
is irreducible in $B$ for each fixed $\bm{\eps}$.
(Note that $\Delta^\textnormal{diag}(P,\bm{c})/3^{\mf{e}_m}$ is a product of such $P$'s.)

But $P$ is \emph{linear} in $E_m$,
with constant term $\pm d_m^3E_1\cdots E_{m-1}$
coprime to the ``slope'' $d_1^3E_2\cdots E_{m-1}\pm\dots\pm d_{m-1}^3E_1\cdots E_{m-2}$.
So by Gauss' lemma, $P$ is irreducible in $B$.
\end{proof}

\begin{remark}
One might try to prove Lemma~\ref{LEM:bi-homogeneous-factored-polynomial-is-irreducible} by specializing $P$ and using Proposition~\ref{PROP:some-general-duality-theory},
but tricky issues might then arise
(since specialization can decrease degrees).

\end{remark}

\section{Establishing codimension-two criteria in general}
\label{SEC:general-codimension-two-work}

In this section, we prove our main general result, Theorem~\ref{THM:main-general-result-informal}, via Theorem~\ref{THM:codimension-2-criteria-for-boundedness-of-E_c} and Proposition~\ref{PROP:some-general-duality-theory}.
First we recall some standard singularity types.
For a finite-type scheme $W$ of pure dimension $n-1\geq 0$ over an algebraically closed field $K$,
a singular $K$-point of $W$ is said to be (of type) \emph{$A_1$} if the completed local ring of $W$ at the point is $\cong K[[z_1,\dots,z_n]]/(z_1^2+\dots+z_n^2)$,
and is said to be \emph{non-$A_1$} otherwise;
and a \emph{non-degenerate double point} of $W$ is a $K$-point at which the completed local ring is $\cong K[[z_1,\dots,z_n]]/(Q)$,
where $Q=Q(z_1,\dots,z_n)$ is a quadratic form over $K$
such that the Hessian matrix of $Q$ is invertible over $K$
(i.e.~the linear forms $\partial Q/\partial z_1, \dots, \partial Q/\partial z_n$ are linearly independent over $K$).

\begin{remark}
\label{RMK:A_1-vs-ordinary-vs-nondegenerate-double-points}
In characteristic $p\neq 2$,
the notions of ``$A_1$ singular $K$-point'' and ``non-degenerate double point'' coincide.
(See e.g.~\cite{poonen2020valuation}*{\S4}.)
\end{remark}

\begin{theorem}
\label{THM:codimension-2-criteria-for-boundedness-of-E_c}
In the Main Setup with $k,m,F,V,V_{\bm{c}}$,
suppose $V$ is smooth, $d\defeq \deg{F}\geq 3$, and $k$ is finite of characteristic $p\neq 2$.
Let $\bm{c}\in k^m\setminus\set{\bm{0}}$, and suppose $V_{\bm{c}}$ is $\abs{E}$-bad (see Definition~\ref{DEFN:error-goodness-for-projective-X/F_q} and Proposition~\ref{PROP:amplification-for-projective-varieties}).
Then the following hold:
\begin{enumerate}
    \item Either $V_{\bm{c}}\times_k \ol{k}$ has $\geq 2$ singular $\ol{k}$-points, or it has a non-$A_1$ singular $\ol{k}$-point.
    
    \item Suppose $p\nmid d(d-1)$.
    Then $\disc(F,\bm{c}) = \partial\disc(F,\bm{c})/\partial c_1 = \dots = \partial\disc(F,\bm{c})/\partial c_m = 0$;
    in other words, $\bm{c}$ is a singular point of the scheme $V_{\Aff^m}(\disc(F,-))_{/k}$.
\end{enumerate}
\end{theorem}

In applications, (2) can be easier to use,
but (1) can also be useful via \cite{poonen2020valuation}*{Theorem~1.1}.

We now turn to the proof.
We first prove Theorem~\ref{THM:codimension-2-criteria-for-boundedness-of-E_c}(1).
For $m\geq 5$, we will apply \cite{lindner2020hypersurfaces}*{Theorem~1.2} (which gives a criterion for non-defectiveness of hypersurfaces of dimension $\geq 2$).
For $m=4$, we will instead use a criterion for integrality of plane curves:

\begin{lemma}
\label{LEM:criterion-for-geometric-integrality-of-plane-curves}
Let $C\belongs \PP^2$ be a hypersurface of degree $\geq 3$ over an algebraically closed field $K$.
Suppose $\abs{C_{\map{sing}}}$ consists of a non-degenerate double point.
Then $C$ is integral.
\end{lemma}

\begin{proof}
Suppose not.
Then $C = V_{\PP^2}(f_1f_2)$, with $f_1,f_2\in K[x,y,z]$ homogeneous of degrees $d_1,d_2\geq 1$.
So $C_{\map{sing}}\contains V(f_1,f_2)$.
Say $C_{\map{sing}}$ is supported at $p\defeq [0:0:1]$.
Then $\abs{V(f_1,f_2)} = \set{p}$, since $V(f_1,f_2)\neq \emptyset$ by dimension theory.
Now let $R\defeq K[[x,y]]$ and $\mf{m}\defeq (x,y)$; then $f_1,f_2\in \mf{m}$.
By our hypothesis on $C_{\map{sing}}$, there exists a continuous $K$-algebra isomorphism $\phi\maps R/(xy)\to R/(f_1f_2)$,
whence $(f_1,f_2) = \mf{m}$ (by a standard calculation).
So by B\'{e}zout, $d_1d_2 = \dim_K(R/(f_1,f_2)) = 1$,
whence $\deg(f_1f_2) = d_1+d_2 = 2$, which is absurd.
\end{proof}

\begin{proof}
[Proof of Theorem~\ref{THM:codimension-2-criteria-for-boundedness-of-E_c}(1)]
For $m=3$, there is nothing to prove, since all $0$-dimensional projective $k$-schemes are $k$-finite (i.e.~finite over $\Spec{k}$), and thus trivially $\abs{E}$-good.
So from now on, suppose $m\geq 4$.
Let $m_\ast\defeq m-3$, so that $V_{\bm{c}}$ is isomorphic to a hypersurface $X\belongs \PP^{1+m_\ast}_k$.
By assumption, $V_{\bm{c}}$ (and thus $X$) is $\abs{E}$-bad, so by the final sentence of Proposition~\ref{PROP:amplification-for-projective-varieties}, there exists $\ell\neq p$ with $\dim H^{1+m_\ast}(X)\neq \dim H^{1+m_\ast}(\PP^{m_\ast}_k)$.
Let $K\defeq \ol{k}$.

Now suppose $m=4$, i.e.~$\dim{X}=m_\ast=1$.
Then $\dim H^2(X)\neq \dim H^2(\PP^1_k) = 1$, so by \cite{poonen2017rational}*{Corollary~7.5.21}, $X_K$ must be reducible.
But $X_K$ is a hypersurface in $\PP^2_K$ of degree $\deg{F}\geq 3$.
So Lemma~\ref{LEM:criterion-for-geometric-integrality-of-plane-curves} and Remark~\ref{RMK:A_1-vs-ordinary-vs-nondegenerate-double-points} give the desired result.

Now suppose $m\geq 5$.
If $X_K$ has $\geq 2$ singular $K$-points or a non-$A_1$ singular $K$-point, then we are done, so suppose otherwise (for contradiction).
Then $p\neq 2$, the hypersurface $X_K$ has \emph{defect} in $\ell$-adic \'{e}tale cohomology, the ambient dimension $1+m_\ast$ is $\geq 3$, and $\abs{(X_K)_{\map{sing}}}$ (if nonempty) consists of isolated $A_1$ singular $K$-points,
so \cite{lindner2020hypersurfaces}*{Theorem~1.2} yields the inequality $2\cdot \#X_{\map{sing}}(K) \geq \deg(X_K) \geq 3$, contradicting $\#X_{\map{sing}}(K)\leq 1$.
\end{proof}

The proof of Theorem~\ref{THM:codimension-2-criteria-for-boundedness-of-E_c}(2) will make use of some general duality theory.
In the Main Setup with $k,m,F,V,V_{\bm{c}}$,
suppose $V$ is smooth, $d\defeq \deg{F}\geq 2$, and $k$ is an arbitrary field of characteristic $p\geq 0$.
When working in the dual projective space $(\PP^{m-1}_k)^\vee$, we will abbreviate $[c_1:\dots:c_m]\in (\PP^{m-1}_k)^\vee(\ol{k})$ as $[\bm{c}]$.
Let $V^\vee$ denote the \emph{dual variety} of $V$, i.e.~the scheme-theoretic image of $V$ under
the \emph{polar map} $[\grad{F}]$ (the rational map $\PP^{m-1}_k\dashrightarrow (\PP^{m-1}_k)^\vee$ sending $[\bm{x}]=[x_1:\dots:x_m]$ to $[\grad{F}(\bm{x})]=[\partial F/\partial x_1:\dots:\partial F/\partial x_m]$;
it is defined on all of $V$, since $V$ is smooth).
Let $\gamma\maps V\to V^\vee$ denote the \emph{Gauss map}, defined by restricting the polar map to $V$.
Let $(V^\vee)_\textnormal{sm}$ denote the largest open subscheme of $V^\vee$ that is smooth over $k$.
Let $I=(\Delta)$ denote the principal ideal of $k[c_1,\dots,c_m]$ generated by the polynomial $\Delta=\disc(F,-)=\disc(F,\bm{c})\in k[c_1,\dots,c_m]$,
and let $\sqrt{I}$ denote the radical of $I$.

\begin{proposition}
\label{PROP:some-general-duality-theory}
Assume the setting above.
The schemes $V,V^\vee$ are geometrically integral hypersurfaces,
the polynomial $\Delta$ is nonconstant,
and $V^\vee$ is defined in $(\PP^{m-1}_k)^\vee$ by the homogeneous ideal $\sqrt{I}$.
If $p\nmid d(d-1)$,
then the following hold:
\begin{itemize}
    \item $\gamma$ is $k$-birational (i.e.~birational over $k$),
    
    \item $(V^\vee)_\textnormal{sm}(\ol{k})
    = \set{[\bm{c}]\in (\PP^{m-1}_k)^\vee(\ol{k}): \abs{((V_{\ol{k}})_{\bm{c}})_{\map{sing}}}\;\textnormal{consists of a non-degenerate double point}}$,
    
    \item the hypersurface $V(\Delta)\belongs (\PP^{m-1}_k)^\vee$ equals $V^\vee$,
    
    \item the polynomial $\Delta$ is irreducible over $\ol{k}$,
    and
    
    \item $\dim((V^\vee)_{\map{sing}}) = -\bm{1}_{d=2} + (m-3)\cdot \bm{1}_{d\geq 3}$.
\end{itemize}
\end{proposition}

\begin{proof}
Since $V_{\ol{k}}$ is smooth of dimension $m-2\geq 1$, it is certainly integral.
In particular, $V_{\ol{k}}$ is reduced, so $V^\vee$ (resp.~$(V_{\ol{k}})^\vee$) is
the space $[\grad{F}](\abs{V})$ (resp.~$[\grad{F}](\abs{V_{\ol{k}}})$), equipped with the reduced induced scheme structure.
Therefore,
$(V^\vee)_{\ol{k}}$ coincides with $(V_{\ol{k}})^\vee$
and is integral,
and
\begin{equation*}
V^\vee(\ol{k})
= \set{[\bm{c}]\in (\PP^{m-1}_k)^\vee(\ol{k}):
((V_{\ol{k}})_{\bm{c}})_{\map{sing}}\neq \emptyset}.
\end{equation*}
So by Proposition-Definition~\ref{PROPDEFN:discriminant-of-hyperplane-sections}, $V(\sqrt{I}) = V^\vee$.
But $\dim V^\vee \leq \dim V$, so $\Delta$ must be nonzero (and thus nonconstant, by degree considerations), whence $V^\vee$ is a hypersurface.

Say $p\nmid d$.
Then $[\grad{F}]$ is defined everywhere, since $V$ is smooth and $d\cdot F = \bm{x}\cdot \grad{F}$.
So $[\grad{F}]\maps \PP^{m-1}_k\to (\PP^{m-1}_k)^\vee$ is a finite surjective morphism (as is $\gamma$).
But $\gamma^\ast\mcal{O}_{V^\vee}(1)\cong \mcal{O}_V(d-1)$,
so the projection formula (along $\gamma$) for intersection numbers then yields $d(d-1)^{m-2} = \deg(\gamma)\deg(V^\vee)$.
Suppose from now on that $p\nmid d(d-1)$; then $\gamma$ is separable (i.e.~the function field of $V$ is separable over the function field of $V^\vee$).
So $\gamma$ is \'{e}tale at some closed point $x\in V$.
The composition $V\to V^\vee\inject (\PP^{m-1}_k)^\vee$ is then unramified at $x$,
so by \cite{deligne1973groupes}*{Katz's Expos\'{e}~XVII, Exemple~3.4, Proposition~3.5, and (3.5.2)--(3.5.3)} (applied to $V_{\ol{k}}\belongs \PP^{m-1}_{\ol{k}}$),
\begin{enumerate}
    \item $(V^\vee)_\textnormal{sm}(\ol{k})$ has the desired alternative description;
    and
    \item if we write $V^\vee = V_{(\PP^{m-1}_k)^\vee}(F^\vee)$, then $\gamma$ and $[\grad{F^\vee}]\maps (\PP^{m-1}_k)^\vee\dashrightarrow \PP^{m-1}_k$ define inverse morphisms between $\gamma^{-1}((V^\vee)_\textnormal{sm})$ and $(V^\vee)_\textnormal{sm}$, after base change to $\ol{k}$.
\end{enumerate}
Since $\ol{k}/k$ is flat, (2) remains true over $k$.
And $(V^\vee)_\textnormal{sm}$ is dense in $V^\vee$, since $(V^\vee)_{\ol{k}}$ is reduced.
So $\dim((V^\vee)_{\map{sing}})\leq \dim(V^\vee)-1 = m-3$, and $\gamma$ is $k$-birational.

In particular, $\deg(\gamma) = 1$, so $\deg(V^\vee) = d(d-1)^{m-2}$.
But $V^\vee = V(\sqrt{I})$, and $\deg(\Delta) = d(d-1)^{m-2}$ by Proposition-Definition~\ref{PROPDEFN:discriminant-of-hyperplane-sections}.
So $V(\Delta) = V^\vee$.
Thus $\Delta$ is irreducible over $\ol{k}$.

It remains to compute $\dim((V^\vee)_{\map{sing}})$.
Let $\bm{M}(\bm{x})$ denote the Hessian matrix of $F$.
If $d=2$, then $\bm{M}(\bm{x})$ is an invertible constant matrix (since $V$ is smooth and $p\nmid d$), so $[\grad{F}]$ is an isomorphism, whence $V^\vee\cong V$ and $\dim((V^\vee)_{\map{sing}}) = \dim(V_{\map{sing}}) = -1$.

Now suppose $d\geq 3$.
Let $[\bm{a}]\in V(\ol{k})$, so that $\bm{a}^t \grad{F}(\bm{a}) = d\cdot F(\bm{a}) = 0$.
Since $\bm{M}(\bm{a}) \bm{a} = (d-1) \grad{F}(\bm{a})$ and $\bm{a}^t \bm{M}(\bm{a}) \bm{v} = (d-1) \grad{F}(\bm{a}) \cdot \bm{v}$,
the formula $\bm{v}\mapsto \bm{M}(\bm{a}) \bm{v}$ defines a $\ol{k}$-linear map $\phi_{\bm{a}}\maps \ol{k}^m/\bm{a}\ol{k} \to \ol{k}^m/\grad{F}(\bm{a})\ol{k}$, whose kernel is orthogonal to $\grad{F}(\bm{a})$.
Here we may identify $\ol{k}^m/\bm{a}\ol{k}$ and $\ol{k}^m/\grad{F}(\bm{a})\ol{k}$ with the tangent spaces $T_{[\bm{a}]} \PP^{m-1}_{\ol{k}}$ and $T_{[\grad{F}(\bm{a})]} (\PP^{m-1}_{\ol{k}})^\vee$, respectively,
and the orthogonal complement of $\grad{F}(\bm{a})$ in $\ol{k}^m/\bm{a}\ol{k}$ with $T_{[\bm{a}]} V_{\PP^{m-1}}(\grad{F}(\bm{a})\cdot \bm{x})_{/\ol{k}}$.
Therefore, the singularity $[\bm{a}]$ of $(V_{\ol{k}})_{[\grad{F}(\bm{a})]}$ is a non-degenerate double point if and only if $\phi_{\bm{a}}$ is injective, which occurs if and only if $\det(\bm{M}(\bm{a}))\neq 0$.
So $\gamma(R)\belongs \abs{(V^\vee)_{\map{sing}}}$, where $R$ denotes the closed subset of $V$ cut out by $\det(\bm{M}(\bm{x})) = 0$.
But $\dim(R)\geq \dim(V)-1 = m-3$, since $\det(\bm{M}(\bm{x}))$ is nonconstant.
So $\dim((V^\vee)_{\map{sing}})\geq m-3$, since $\gamma$ is finite.
Thus $\dim((V^\vee)_{\map{sing}}) = m-3$.
\end{proof}

\begin{proof}
[Proof of Theorem~\ref{THM:codimension-2-criteria-for-boundedness-of-E_c}(2)]
By Theorem~\ref{THM:codimension-2-criteria-for-boundedness-of-E_c}(1) and Remark~\ref{RMK:A_1-vs-ordinary-vs-nondegenerate-double-points}, $\abs{((V_{\ol{k}})_{\bm{c}})_{\map{sing}}}$ is \emph{not} a singleton consisting of a non-degenerate double point.
Now say $p\nmid d(d-1)$.
Then $[\bm{c}]\in (V^\vee)_{\map{sing}}$, by Proposition~\ref{PROP:some-general-duality-theory}.
But also by Proposition~\ref{PROP:some-general-duality-theory}, $V^\vee$ is the zero locus of a polynomial $F^\vee\mid \Delta$.
(In fact, $\Delta \in k^\times\cdot F^\vee$, since $p\nmid d(d-1)$.)
So $\Delta$ vanishes, and is singular, at $\bm{c}$.
\end{proof}

\begin{proof}
[Proof of Theorem~\ref{THM:main-general-result-informal}]
Immediate from Theorem~\ref{THM:codimension-2-criteria-for-boundedness-of-E_c}(2) and Proposition~\ref{PROP:some-general-duality-theory}.
\end{proof}

\section{Establishing near-optimal criteria in special cases}
\label{SEC:establishing-near-dichotomy-for-quadrics-and-low-dimensional-cubics}

Call a projective complete intersection of multi-degree $\bm{d}$ a \emph{$\bm{d}$-CI}.
(In the Main Setup, for instance, $V_{\bm{c}}$ is a $(d, 1)$-CI in $\PP^{m-1}_k$ if and only if $\dim(V_{\bm{c}}) = m-3$.)
In this section, we prove our main specialized result, Theorem~\ref{THM:main-special-result-informal}, via Theorem~\ref{THM:near-dichotomy-for-cubic-curves-and-threefolds}.
As preparation,
we first briefly discuss the essentially classical situation of a $(2)$-CI, i.e.~a projective quadric.

\begin{definition}
Let $k$ be a base field.
An embedded projective $k$-scheme is a \emph{cone} if and only if it is ``missing a variable'' after some $k$-linear change of coordinates.
(E.g.~$V_{\PP^1}((x_1-x_2)^3)_{/k}$ is a cone.)
A \emph{cone over} a projective $k$-scheme $Y$ is a projective cone over $Y$, with vertex a $k$-point.
An \emph{iterated cone} is obtained by taking cones one or more times.
\end{definition}

\begin{proposition}
\label{PROP:dichotomy-for-quadrics}
Let $X$ be a $(2)$-CI of dimension $\geq 0$ over a finite field $k$.
Then $X$ is stably $\abs{E}$-bad
\emph{if and only if} $X_{\ol{k}}$ is an iterated cone over a smooth $(2)$-CI of dimension $d\in [0, \dim(X)-1]$ with $2\mid d$;
\emph{if and only if} $X$ is singular and $\dim(X_{\map{sing}})\equiv \dim(X)-1\bmod{2}$.
\end{proposition}

\begin{remark}
If $2\nmid \card{k}$, then the adverb ``stably'' can be removed, and the cone condition is equivalent to ``$X$ is of the form $V(Q)$ with $2\mid \rank(Q)\in [2, \dim(X)+1]$''.
(Here $\rank(Q)$ denotes the rank of the Hessian matrix of $Q$.)
\end{remark}

Proposition~\ref{PROP:dichotomy-for-quadrics} (when combined with Theorem~\ref{THM:Zak's-principle}, Proposition~\ref{PROP:amplification-for-projective-varieties}, and \cite{debarre2003lines}*{Lemma~3}, say) implies that in the Main Setup, if $k$ is finite, $V$ is smooth, $\deg{F}=2$, and $2\mid m$, then $V_{\bm{c}}$ is $\abs{E}$-bad if and only if $\disc(F,\bm{c})=0$;
if and only if $V_{\bm{c}}\times_k \ol{k}$ contains an $\frac{m-2}{2}$-plane in $\PP^{m-1}_{\ol{k}}$.
In particular, the assumption $\deg{F}\geq 3$ in Theorems~\ref{THM:main-general-result-informal} and~\ref{THM:codimension-2-criteria-for-boundedness-of-E_c} is essential when $2\mid m$.

To state Theorem~\ref{THM:near-dichotomy-for-cubic-curves-and-threefolds}, we need a definition:

\begin{definition}
\label{DEFN:cubic-scroll-of-dimension-d}
Let $d\geq 1$ be an integer.
Over an algebraically closed field, a \emph{cubic $d$-scroll} is
an embedded projective scheme $\Sigma$, integral of dimension $d$ and degree $3$, with $\dim(\Span(\Sigma)) = d+2$.
(There are other equivalent definitions; see e.g.~\cite{eisenbud1987varieties}.)
Here $\Span(\Sigma)$ denotes the \emph{projective span} of $\Sigma$, i.e.~the smallest projective space containing $\Sigma$.
\end{definition}

\begin{theorem}
\label{THM:near-dichotomy-for-cubic-curves-and-threefolds}
Let $n\in \set{2,4}$.
Let $X$ be a $(3)$-CI in $\PP^n$ over a finite field $k$.
Suppose $X_{\ol{k}}$ is not a cone over a cone.
If $n=2$, then $X$ is $\abs{E}$-bad if and only if $X_{\ol{k}}$ contains a line in $\PP^2_{\ol{k}}$.
Now suppose $n=4$, and consider the following four conditions:
\begin{enumerate}
    \item\label{ITEM:cubic-threefold-stably-E-bad} $X$ is stably $E$-bad (see Definition~\ref{DEFN:error-goodness-for-projective-X/F_q} and Proposition~\ref{PROP:amplification-for-projective-varieties});
    \item\label{ITEM:cubic-threefold-stably-|E|-bad} $X$ is stably $\abs{E}$-bad;
    \item\label{ITEM:cubic-threefold-potentially-contains-a-special-subvariety} $X_{\ol{k}}$ contains either a plane, or a singular cubic $2$-scroll, in $\PP^4_{\ol{k}}$;
    and
    \item\label{ITEM:cubic-threefold-potentially-contains-a-nonempty-open-subscheme-of-a-special-CI} there exists a $(2,2)$-CI of the form $V(Q_1,Q_2)\belongs \PP^4_{\ol{k}}$, with $V(Q_1)_{\map{sing}}\cap V(Q_2)_{\map{sing}}\neq \emptyset$, such that $X_{\ol{k}}$ contains a nonempty open subscheme of $V(Q_1,Q_2)$.
\end{enumerate}
In general,
(\ref{ITEM:cubic-threefold-stably-E-bad})--(\ref{ITEM:cubic-threefold-stably-|E|-bad}) are equivalent,
(\ref{ITEM:cubic-threefold-stably-|E|-bad}) implies (\ref{ITEM:cubic-threefold-potentially-contains-a-special-subvariety}),
and (\ref{ITEM:cubic-threefold-potentially-contains-a-special-subvariety})--(\ref{ITEM:cubic-threefold-potentially-contains-a-nonempty-open-subscheme-of-a-special-CI}) are equivalent.
If $\dim(X_{\map{sing}})\leq 0$,
then (\ref{ITEM:cubic-threefold-stably-E-bad})--(\ref{ITEM:cubic-threefold-potentially-contains-a-nonempty-open-subscheme-of-a-special-CI}) are equivalent.
\end{theorem}

\begin{remark}
\label{RMK:analytic-auxiliary-polynomial-complexity-direction-of-inquiry}
In (\ref{ITEM:cubic-threefold-potentially-contains-a-nonempty-open-subscheme-of-a-special-CI}), and in similar places below, we have opted to use somewhat analytic language in lieu of more typical geometric language:
(\ref{ITEM:cubic-threefold-potentially-contains-a-nonempty-open-subscheme-of-a-special-CI}) says that the hypersurface $X_{\ol{k}}=V(C)_{/\ol{k}}$ (where $C$ is a cubic form defining $X$) ``almost'' contains $V(Q_1,Q_2)$, in the sense that there exists a polynomial $A\notin \sqrt{(Q_1,Q_2)}$ with $A\cdot C \in (Q_1,Q_2)$.
We include this formulation because it may suggest a new line of inquiry, where instead of identifying explicit structure one tries to produce auxiliary polynomial conditions involving $X$.
In the specific present setting however,
one could replace the phrase ``nonempty open subscheme'' with the more usual ``irreducible component (with reduced induced scheme structure)'',
and also replace ``$V(Q_1)_{\map{sing}}\cap V(Q_2)_{\map{sing}}\neq \emptyset$'' with ``$V(Q_1,Q_2)$ is a cone''.
It is also worth noting that if $X_{\ol{k}}$ contains a plane, then it must completely contain a $(2,2)$-CI,
but not all $(2,2)$-CI's on $X_{\ol{k}}$ need be cones.
\end{remark}

\begin{remark}
The equivalence of (\ref{ITEM:cubic-threefold-stably-E-bad})--(\ref{ITEM:cubic-threefold-stably-|E|-bad}) suggests that (stable) $\abs{E}$-badness might be explained by ``excess'' points from ``special'' subvarieties, which we have tried to pinpoint in (\ref{ITEM:cubic-threefold-potentially-contains-a-special-subvariety})--(\ref{ITEM:cubic-threefold-potentially-contains-a-nonempty-open-subscheme-of-a-special-CI}).
Theorem~\ref{THM:near-dichotomy-for-cubic-curves-and-threefolds} is close to a complete dichotomy.
To ``complete'' it (for $n=4$), one would need to analyze the case $\dim(X_{\map{sing}})\geq 1$, which might be tricky (see e.g.~the proof of Lemma~\ref{LEM:associated-(2,3)-CI-controls-potential-goodness-of-singular-cubic-threefold}).
\end{remark}

\begin{remark}
The case when $X_{\ol{k}}$ is a cone over a cone can still be fully analyzed,
but it is less interesting than the opposite case.
Our methods can also be used to show that a $(3)$-CI in $\PP^3_k$ is stably $\abs{E}$-bad if and only if $X_{\ol{k}}$ is either reducible, or a cone over a smooth plane cubic curve;
we omit this from Theorem~\ref{THM:near-dichotomy-for-cubic-curves-and-threefolds} because it has a different flavor.
\end{remark}

Theorem~\ref{THM:main-special-result-informal} readily follows from Theorem~\ref{THM:Zak's-principle}, Proposition~\ref{PROP:amplification-for-projective-varieties}, and Theorem~\ref{THM:near-dichotomy-for-cubic-curves-and-threefolds}.
So we devote the rest of \S\ref{SEC:establishing-near-dichotomy-for-quadrics-and-low-dimensional-cubics} to the proof of Proposition~\ref{PROP:dichotomy-for-quadrics} and Theorem~\ref{THM:near-dichotomy-for-cubic-curves-and-threefolds}.
Proposition~\ref{PROP:dichotomy-for-quadrics} stems from the following general fact (which we state without proof):

\begin{proposition}
[Routine; see e.g.~\cite{wang2022thesis}*{Proposition~5.1.23}]
\label{PROP:cone-construction-multiplies-E-by-q}
Let $k$ be a finite field.
If $C(Y)$ is a cone over a projective $k$-scheme $Y$, then $E(C(Y)) = \card{k}\cdot E(Y)$.
\end{proposition}

\begin{proof}
[Proof sketch for Proposition~\ref{PROP:dichotomy-for-quadrics}]
The key ingredients are the following:
\begin{enumerate}
    \item If $X$ is smooth of dimension $d\geq 0$, then $E(X) = 0$ if $2\nmid d$, and $E(X) = \pm \card{k}^{d/2}$ if $2\mid d$.
    (One can prove this using either the Weil conjectures
    or elementary tools.)
    
    \item If $X_{\ol{k}}$ is non-reduced, then $E(X)=0$ and $\dim(X_{\map{sing}}) = \dim(X)$.
    If $X_{\ol{k}}$ is reduced, then $X_{\ol{k}}$ is singular if and only if it is an iterated cone over a smooth $(2)$-CI of dimension $d\geq 0$; and in this case, we must have $\dim(X_{\map{sing}}) = \dim(X)-d-1$.
\end{enumerate}
By (1)--(2) and Proposition~\ref{PROP:cone-construction-multiplies-E-by-q}, $X$ is potentially $\abs{E}$-good if and only if $X_{\ol{k}}$ is either smooth or non-reduced, or an iterated cone over a smooth $(2)$-CI of dimension $d\in [0, \dim(X)-1]$ with $2\nmid d$.
The result follows from (2).
\end{proof}

Our proof of Theorem~\ref{THM:near-dichotomy-for-cubic-curves-and-threefolds} uses (for $n=4$) base change and some classical geometry, including some classification results over $\ol{k}$.
(Base change is a priori ``legal'' by Lemma~\ref{LEM:potential-and-stable-logic}.)
One basic idea in the proof is that a reducible curve of degree $2\cdot 3 = 6$ must have an irreducible component of degree $\le 3$; however, the full proof is rather more technical.
Originally, we sought to use an extension of \cite{bombieri1967local}'s ``conic bundle'' method;
see \cite{wang2022thesis}*{\S5.3, especially Remark~5.3.16} for details.
This approach inspired condition~(\ref{ITEM:cubic-threefold-potentially-contains-a-nonempty-open-subscheme-of-a-special-CI}), in fact in an ``effective'' form.
But our present approach is overall more efficient in the singular case.

It would be interesting to find different proofs of Theorem~\ref{THM:near-dichotomy-for-cubic-curves-and-threefolds} that work directly over $k$, or that generalize naturally.
Condition~(\ref{ITEM:cubic-threefold-potentially-contains-a-nonempty-open-subscheme-of-a-special-CI}) in Theorem~\ref{THM:near-dichotomy-for-cubic-curves-and-threefolds} seems especially suggestive as to what one might try more generally
(e.g.~if one were to try using the method of auxiliary polynomials).
One might also try using mixed Hodge theory, in the spirit of e.g.~\cites{dimca1990betti,kloosterman2022maximal}.

In any case, we now begin with a classical ``rationality''-type idea (cf.~\cite{dolgachev2016corrado}*{\S1}).
\begin{proposition}
\label{PROP:using-k-rational-singular-point-to-count-points}
Let $X\belongs \PP^n$ be a $(3)$-CI over a finite field $k$, where $n\geq 2$.
Assume $[0:\dots:0:1]\in X_{\map{sing}}$.
Then $X = V_{\PP^n}(f_2x_{n+1}+f_3)$ for some $f_2, f_3\in k[x_1,\dots,x_n]$, homogeneous of respective degrees $2, 3$, with $(f_2,f_3)\neq (0,0)$.
Furthermore,
\begin{equation}
\label{EQN:rationality-based-equation-for-E(X)}
E(X) = \card{k}\cdot E(V_{\PP^{n-1}}(f_2,f_3)) - E(V_{\PP^{n-1}}(f_2)) + \card{k}^{n-1}\cdot \bm{1}_{\dim V_{\PP^{n-1}}(f_2,f_3) = \dim V_{\PP^{n-1}}(f_2)}.
\end{equation}
(As in Definition~\ref{DEFN:analytic-and-combinatorial-notation}, we let $\bm{1}_A$ denote the \emph{indicator value} of an event $A$.)
\end{proposition}

\begin{proof}
The first part is clear.
So there are two kinds of points $[\bm{x}]\in X(k)$: (i) those with $f_2\neq 0$ and $x_{n+1} = -f_3/f_2$, and (ii) those with $f_2=0$ and $f_3=0$.
Therefore
\begin{equation*}
\card{X(k)} = \card{(\PP^{n-1}\setminus V_{\PP^{n-1}}(f_2))(k)} + \card{C(V_{\PP^{n-1}}(f_2,f_3))(k)},
\end{equation*}
where $C(-)$ denotes a cone.
Proposition~\ref{PROP:cone-construction-multiplies-E-by-q}, and casework on the quantity $\dim V_{\PP^{n-1}}(f_2)-\dim V_{\PP^{n-1}}(f_2,f_3)\in \set{0,1}$, then lead to the desired equality.
\end{proof}

To study $E(X)$ using Proposition~\ref{PROP:using-k-rational-singular-point-to-count-points}, we will need to analyze low-degree complete intersections in some detail.
We will frequently (and sometimes implicitly) use the fact that projective complete intersections are equidimensional
(and in fact Cohen--Macaulay;
see e.g.~\cite{stacks-project}*{\href{https://stacks.math.columbia.edu/tag/00SB}{Tag~00SB} and \href{https://stacks.math.columbia.edu/tag/00N9}{Tag~00N9}}),
as well as the following standard technical result:

\begin{proposition}
\label{PROP:locally-Noetherian-generically-reduced-and-CM-implies-reduced}
Let $X$ be a Noetherian, Cohen--Macaulay scheme.
Let $U\belongs X$ be an open subscheme that is topologically dense.
Suppose $U$ is regular or reduced.
Then $X$ is reduced.
\end{proposition}

\begin{proof}
[Proof sketch]
Any regular local ring is an integral domain (and in particular, reduced), so we may assume $U$ is reduced.
Now if $\eta$ is a generic point of an irreducible component of $U$, then the $0$-dimensional local ring $\mcal{O}_{U, \eta}$ is reduced, and thus a field (i.e.~regular) \cite{stacks-project}*{\href{https://stacks.math.columbia.edu/tag/00EU}{Tag~00EU}}.
Since $\abs{U}$ is dense in $\abs{X}$,
it follows that $X$ is regular in codimension $0$ (in the sense of \cite{stacks-project}*{\href{https://stacks.math.columbia.edu/tag/033Q}{Tag~033Q}}).
By \cite{stacks-project}*{\href{https://stacks.math.columbia.edu/tag/0342}{Tag~0342} and \href{https://stacks.math.columbia.edu/tag/0344}{Tag~0344}}, we conclude that $X$ is reduced.
\end{proof}

We will also repeatedly use the following lemma describing the low-degree irreducible components of non-integral $(2,2)$-CI's and $(2,3)$-CI's:

\begin{lemma}
\label{LEM:description-of-low-degree-components-of-a-CI}
Let $n\geq 2$.
Let $K$ be an algebraically closed field.
Let $Y$ be a $(d,e)$-CI of the form $V(A,B)\belongs \PP^n_K$, with $\deg{A}=d=2$ and $\deg{B}=e\in \set{2,3}$.
Then $Y$ is non-integral if and only if it has an irreducible component of degree $\leq 3$.
Let $Z$ be any such component, equipped with the reduced induced scheme structure.
Then the following dichotomy holds:
\begin{enumerate}
    \item If $\deg{Z}\leq 2$, or $e=3$ and $A$ is reducible, then $\dim(\Span(Z))\leq n-1$.
    
    \item If $\deg{Z}=3$, and $e=2$ or $A$ is irreducible, then $\dim(\Span(Z)) = n$.
\end{enumerate}
Furthermore, $\dim(\Span(Z))\leq n-1$ if and only if $Z$ is a $(1,\deg{Z})$-CI in $\PP^n_K$.
\end{lemma}

\begin{proof}
The first part is clear, since $\deg{Y} = de\leq 6$.
Now fix $Z$.
Since $\dim{Z} = n-2$, and $Z$ is integral, the final sentence is clear: both conditions are equivalent to ``$Z$ lies in an $(n-1)$-plane in $\PP^n_K$ (scheme-theoretically)''.
So it remains to prove (1)--(2).
If $\deg{Z}\leq 2$, then $\dim(\Span(Z))\leq \dim(Z)+1 = n-1$ by \cite{eisenbud1987varieties}*{Proposition~0},
so (1)--(2) hold here.

Now suppose $\deg{Z}=3$.
If $A$ is reducible, then since $Z$ is integral, there must exist a linear factor $L\mid A$ such that $V(L,B)\contains Z$; and thus $e=3$ and $Z = V(L,B)\belongs V(L)$.
Conversely, suppose $Z = V(L,C)$ for some nonzero linear form $L$ and cubic form $C$ (with $L\nmid C$).
Then $V(A,B)\contains V(L,C)$.
A short calculation in $V(L)\cong \PP^{n-1}_K$ then shows that $L\mid A$ and $e=3$ (or else $L\mid A,B$, which is impossible).

Thus we have shown that if $\deg{Z}=3$, then $Z$ is a $(1,3)$-CI if and only if $A$ is reducible, in which case $e=3$ must hold.
This proves (1)--(2) when $\deg{Z}=3$.
\end{proof}

We now prove a technical result leading up to the key Lemma~\ref{LEM:associated-(2,3)-CI-controls-potential-goodness-of-singular-cubic-threefold}.
\begin{lemma}
\label{LEM:lower-bound-on-singular-locus-in-terms-of-any-associated-(2,3)-CI}
In the setting of Proposition~\ref{PROP:using-k-rational-singular-point-to-count-points}, suppose $Y\defeq V_{\PP^{n-1}}(f_2,f_3)$ is a $(2,3)$-CI.
Let $Q\defeq V_{\PP^{n-1}}(f_2)$.
Then the following hold:
\begin{enumerate}
    \item If $\dim(Q_{\map{sing}})\leq \dim(Y_{\map{sing}})-1$, then $\dim(X_{\map{sing}})\geq \dim(Y_{\map{sing}})$.
    
    \item If $x_n\mid f_2$, then $\dim(X_{\map{sing}})\geq \dim(V_{\PP^{n-1}}(x_n,f_3)_{\map{sing}})$.
\end{enumerate}
\end{lemma}

\begin{proof}
By the definitions, $X_{\map{sing}} = V_{\PP^n}(f_2,f_3,x_{n+1}\grad(f_2) + \grad(f_3))$.
So the formula $[x_1:\dots:x_n]\mapsto [x_1:\dots:x_n:-\grad(f_3)/\grad(f_2)]$ defines an embedding $Y_{\map{sing}}\setminus Q_{\map{sing}}\inject X_{\map{sing}}$.
Claim~(1) follows by dimension theory.
On the other hand, if $x_n\mid f_2$, then by a short calculation,
\begin{equation*}
X_{\map{sing}}\contains V_{\PP^n}(x_n, f_3, \grad(f_3\vert_{x_n=0}), x_{n+1}(\partial f_2 / \partial x_n)+(\partial f_3 / \partial x_n)),
\end{equation*}
whence $X_{\map{sing}}$ contains
the intersection of $V_{\PP^n}(x_{n+1}(\partial f_2 / \partial x_n)+(\partial f_3 / \partial x_n))$ with
\begin{equation*}
V_{\PP^n}(x_n, f_3, \grad(f_3\vert_{x_n=0}))
= V_{\PP^n}(x_n,f_3)_{\map{sing}}
= C_{[0:\dots:0:1]}(V_{\PP^{n-1}}(x_n,f_3)_{\map{sing}}).
\end{equation*}
(Here $\grad(f_3\vert_{x_n=0})$ denotes the gradient of $f_3\vert_{x_n=0}$ with respect to $x_1,\dots,x_{n-1}$,
and $C_v(Z)$ denotes the cone over $Z$ with vertex $v$.)
Claim~(2) then follows by dimension theory.
\end{proof}

The following lemma represents significant progress towards Theorem~\ref{THM:near-dichotomy-for-cubic-curves-and-threefolds}:

\begin{lemma}
\label{LEM:associated-(2,3)-CI-controls-potential-goodness-of-singular-cubic-threefold}
Let $X\belongs \PP^4$ be a $(3)$-CI over a finite field $k$.
Suppose $X_{\ol{k}}$ is not a cone over a cone.
Assume $[0:0:0:0:1]\in X_{\map{sing}}$, and let $f_2,f_3$ be as in Proposition~\ref{PROP:using-k-rational-singular-point-to-count-points}.
Assume $f_2,f_3$ are nonzero and coprime.
Then among the following conditions, (1)--(2) are equivalent, and (2) implies (3):
\begin{enumerate}
    \item $X$ is stably $E$-bad;
    \item $X$ is stably $\abs{E}$-bad;
    and
    \item $V_{\PP^3}(f_2,f_3)_{/\ol{k}}$ contains a $(1,1)$-CI, a $(1,2)$-CI, or a cubic $1$-scroll, in $\PP^3_{\ol{k}}$.
\end{enumerate}
Furthermore, if $\dim(X_{\map{sing}})\leq 0$, then (1)--(3) are equivalent.
\end{lemma}

\begin{proof}
Our assumptions on $f_2,f_3$ imply that $Y\defeq V_{\PP^3}(f_2,f_3)$ is a $(2,3)$-CI curve.
In particular, $\dim V_{\PP^3}(f_2) = 1+\dim V_{\PP^3}(f_2,f_3)$,
so for each finite extension $k'/k$, Proposition~\ref{PROP:using-k-rational-singular-point-to-count-points} implies
\begin{equation}
\label{EQN:E-formula-given-k-rational-singular-point-in-(2,3)-CI-case}
E(X_{k'}) = \card{k'}\cdot E(V_{\PP^3}(f_2,f_3)_{/k'}) - E(V_{\PP^3}(f_2)_{/k'}).
\end{equation}

\emph{Case~1: $f_2$ is irreducible over $\ol{k}$.}
Then $E(V_{\PP^3}(f_2)_{/k'}) = O(\card{k'}^{3/2})$ by Lang--Weil.
So if $A\in \set{E, \abs{E}}$, then by eq.~\eqref{EQN:E-formula-given-k-rational-singular-point-in-(2,3)-CI-case},
$X$ is stably $A$-bad if and only if $Y$ is stably $A$-bad; or equivalently (by Lang--Weil), $Y_{\ol{k}}$ is reducible;
and in this situation, $Y_{\ol{k}}$ is non-integral, so condition~(3) holds (by Lemma~\ref{LEM:description-of-low-degree-components-of-a-CI}).
Furthermore, if $\dim(X_{\map{sing}})\leq 0$, then $Y_{\ol{k}}$ is reduced by Proposition~\ref{PROP:locally-Noetherian-generically-reduced-and-CM-implies-reduced} and Lemma~\ref{LEM:lower-bound-on-singular-locus-in-terms-of-any-associated-(2,3)-CI}(1) (since $\dim(V_{\PP^3}(f_2)_{\map{sing}})\leq 0$), so one can reverse the implications (by Lemma~\ref{LEM:description-of-low-degree-components-of-a-CI}).

\emph{Case~2: $f_2$ is reducible over $\ol{k}$.}
Then by base change (via Lemma~\ref{LEM:potential-and-stable-logic}) and a linear change of coordinates, we may assume that $f_2\in \set{x_4^2, x_3x_4}$.

\emph{Subcase~2.1: $f_2 = x_4^2$.}
Then $x_4\nmid f_3$ (i.e.~$V_{\PP^3}(x_4,f_3)$ is a curve), and eq.~\eqref{EQN:E-formula-given-k-rational-singular-point-in-(2,3)-CI-case} simplifies to $E(X_{k'}) = \card{k'}\cdot E(V_{\PP^3}(x_4,f_3)_{/k'})$.
So if $A\in \set{E, \abs{E}}$,
then $X$ is stably $A$-bad if and only if $V_{\PP^3}(x_4,f_3)$ is stably $A$-bad;
or equivalently (by Lang--Weil), $V_{\PP^3}(x_4,f_3)_{/\ol{k}}$ is reducible;
and in this situation, $Y_{\ol{k}}$ contains a $(1,1)$-CI in $\PP^3_{\ol{k}}$, so condition~(3) holds.
Furthermore, if $\dim(X_{\map{sing}})\leq 0$, then by Proposition~\ref{PROP:locally-Noetherian-generically-reduced-and-CM-implies-reduced} and Lemma~\ref{LEM:lower-bound-on-singular-locus-in-terms-of-any-associated-(2,3)-CI}(2), $V_{\PP^3}(x_4,f_3)_{/\ol{k}}$ is reduced, so one can reverse the implications: if condition~(3) holds, then for degree reasons, $V_{\PP^3}(x_4,f_3)_{/\ol{k}}$ must be non-integral (and thus reducible), since it is not a cubic $1$-scroll.

\emph{Subcase~2.2: $f_2 = x_3x_4$.}
Then $x_3\nmid f_3$ and $x_4\nmid f_3$, and eq.~\eqref{EQN:E-formula-given-k-rational-singular-point-in-(2,3)-CI-case} leads, upon a short calculation, to
\begin{equation}
\label{EQN:E-formula-given-k-rational-singular-point-in-(2,3)-CI-case-with-f_2-split}
\begin{split}
E(X_{k'}) = 1 &- (\card{k'}\cdot \card{V_{\PP^3}(x_3,x_4,f_3)(k')} - \card{V_{\PP^3}(x_3,x_4)(k')}) \\
&+ \sum_{i\in \set{3,4}}(\card{k'}\cdot \card{V_{\PP^3}(x_i,f_3)(k')} - \card{V_{\PP^3}(x_i)(k')}).
\end{split}
\end{equation}

\emph{Subsubcase~2.2.1: $f_3(x_1,x_2,0,0) = 0$.}
Then $f_3\in (x_3,x_4)$, so eq.~\eqref{EQN:E-formula-given-k-rational-singular-point-in-(2,3)-CI-case-with-f_2-split} simplifies to
\begin{equation*}
E(X_{k'}) = O(\card{k'}) - \card{k'}^2 + \card{k'}\cdot \sum_{i\in \set{3,4}} \card{(V_{\PP^3}(x_i,f_3)\setminus V_{\PP^3}(x_3,x_4))(k')}.
\end{equation*}
But since $X = V_{\PP^4}(x_3x_4x_5+f_3)$ is (by assumption) not a cone over a cone, $f_3\notin (x_3x_4) + (x_3,x_4)^3 = (x_3x_4,x_3^3,x_4^3) = (x_3,x_4^3)\cap (x_3^3,x_4)$.
So $f_3\in (x_3,x_4)\setminus (x_i,x_3^3,x_4^3)$ for some $i\in \set{3,4}$, whence $V_{\PP^3}(x_i,f_3)$ is reducible.
By Lang--Weil, it follows that $X$ is potentially $(-E)$-good.
So conditions~(1)--(2) are equivalent.
And $Y\contains V_{\PP^3}(x_3,x_4)$, since $f_3\in (x_3,x_4)$; so condition~(3) always holds.
Furthermore, if $\dim(X_{\map{sing}})\leq 0$, then by Proposition~\ref{PROP:locally-Noetherian-generically-reduced-and-CM-implies-reduced} and Lemma~\ref{LEM:lower-bound-on-singular-locus-in-terms-of-any-associated-(2,3)-CI}(2), each $V_{\PP^3}(x_i,f_3)_{/\ol{k}}$ is reduced, and therefore each has an irreducible component distinct from $\abs{V_{\PP^3}(x_3,x_4)}$; so conditions~(1)--(3) all hold.

\emph{Subsubcase~2.2.2: $f_3(x_1,x_2,0,0) \neq 0$.}
Then eq.~\eqref{EQN:E-formula-given-k-rational-singular-point-in-(2,3)-CI-case-with-f_2-split} simplifies to
\begin{equation*}
E(X_{k'}) = O(\card{k'}) + \sum_{i\in \set{3,4}}\card{k'}\cdot E(V_{\PP^3}(x_i,f_3)_{/k'}).
\end{equation*}
So if $A\in \set{E,\abs{E}}$, then by Lang--Weil,
$X$ is stably $A$-bad if and only if $V_{\PP^3}(x_i,f_3)_{/\ol{k}}$ is reducible for some $i\in \set{3,4}$;
and in this situation, $Y_{\ol{k}}$ contains a $(1,1)$-CI in $\PP^3_{\ol{k}}$, so condition~(3) holds.
Furthermore, if $\dim(X_{\map{sing}})\leq 0$, then each $V_{\PP^3}(x_i,f_3)_{/\ol{k}}$ is reduced by Proposition~\ref{PROP:locally-Noetherian-generically-reduced-and-CM-implies-reduced} and Lemma~\ref{LEM:lower-bound-on-singular-locus-in-terms-of-any-associated-(2,3)-CI}(2), so one can reverse the implications: if condition~(3) holds, then for degree reasons, some $V_{\PP^3}(x_i,f_3)_{/\ol{k}}$ must be non-integral (and thus reducible), since neither of them is a cubic $1$-scroll.

In each case, conditions~(1)--(2) are equivalent, and they imply condition~(3); and if $\dim(X_{\map{sing}})\leq 0$, then conditions~(1)--(3) are equivalent.
\end{proof}

Besides Lemma~\ref{LEM:associated-(2,3)-CI-controls-potential-goodness-of-singular-cubic-threefold}, the following four structural results are also important.

\begin{lemma}
[Cf.~\cite{carlini2008complete}*{Remark~4.2}]
\label{LEM:convert-(1,2)-CI-or-(2,2)-CI-containment-to-(1,1)-CI-containment}
Let $n\geq 2$.
Let $X\belongs \PP^n$ be a $(3)$-CI over a field $k$.
Say $X$ contains a $\bm{d}$-CI in $\PP^n$ with $\bm{d}\in \set{(1,2), (2,2)}$.
Then $X$ contains a $(1,1)$-CI in $\PP^n$.
\end{lemma}

\begin{proof}
Let the $\bm{d}$-CI be of the form $V(P_1,P_2)$, with $\deg{P_1}\in \set{1,2}$ and $\deg{P_2} = 2$.
The homogeneous ideal $(P_1,P_2)$ is saturated, so $X$ is of the form $V(A_1P_1+A_2P_2)$, for some $A_i$ with $\deg{A_i} = 3-\deg{P_i}$.
But then $X$ contains $V(A_1,A_2)$ and $V(P_1,A_2)$, one of which is either a $(1,1)$-CI or a $(1)$-CI.
This suffices.
\end{proof}

\begin{proposition}
\label{PROP:basic-cubic-scroll-facts}
Let $d\geq 1$.
Let $\Sigma$ be a cubic $d$-scroll over an algebraically closed field.
\begin{enumerate}
    \item\label{ITEM:cubic-scroll-contains-a-nonempty-open-subscheme-of-a-(2,2)-CI} $\Sigma$ contains a nonempty open subscheme of a $(2,2)$-CI in $\Span(\Sigma)$.
    
    \item\label{ITEM:singular-cubic-scroll-equivalences} $\Sigma$ is singular if and only if it is a cone, in which case it is a cone over a cubic $(d-1)$-scroll (and we have $d\geq 2$).
\end{enumerate}
\end{proposition}

\begin{proof}
[Proof sketch]
Compare Definition~\ref{DEFN:cubic-scroll-of-dimension-d} and \cite{eisenbud1987varieties}*{definition of a ``variety of minimal degree''}, and then inspect \cite{eisenbud1987varieties}*{Theorem~1, and pp.~4--6 of \S1}.
Note that if $\Sigma$ is a cone, say over $Y$,
then $Y\belongs \Sigma$ is integral of dimension $d-1$ with $\deg{Y} = \deg{\Sigma}$ and $\dim(\Span(Y)) = \dim(\Span(\Sigma)) - 1$, so $Y$ is a cubic $(d-1)$-scroll.
\end{proof}

\begin{lemma}
\label{LEM:any-special-subvariety-is-a-cone-with-vertex-singular}
Let $X\belongs \PP^4$ be a $(3)$-CI over an algebraically closed field $K$.
Suppose $X$ contains a surface $S\belongs \PP^4$ that is either a plane or a singular cubic $2$-scroll.
Then there exists a $K$-point $v\in X_{\map{sing}}$ such that $S$ is a cone with vertex $v$.
\end{lemma}

\begin{proof}
If $S$ is a plane, then it must intersect $X_{\map{sing}}$ (see e.g.~\cite{debarre2003lines}*{Lemma~3}); any $K$-point $v\in S\cap (X_{\map{sing}})$ then suffices.
Now suppose $S$ is a singular cubic $2$-scroll.
By Proposition~\ref{PROP:basic-cubic-scroll-facts}(\ref{ITEM:singular-cubic-scroll-equivalences}), $S$ is a cone.
We may then assume (by a linear change of coordinates) that $S$ is cut out by equations involving only $x_1,\dots,x_4$ (not $x_5$);
geometrically, this means $S$ is a cone with vertex $v = [0:0:0:0:1]$.
Now write $X = V_{\PP^4}(f_0x_5^3+f_1x_5^2+f_2x_5+f_3)$, where the $f_i$ are homogeneous in $x_1,\dots,x_4$ of degree $i$.
Then $S\belongs X$ forces each $f_i$ to lie in the saturated homogeneous ideal defining $S$.
So $S\belongs V_{\PP^4}(f_0,f_1,f_2,f_3)$.
Since $S\neq \emptyset$ and $\dim(\Span(S)) = 4$, it follows that $f_0=f_1=0$.
Thus $X$ is singular at $v$, so $v$ suffices.
\end{proof}

Lemma~\ref{LEM:any-special-subvariety-is-a-cone-with-vertex-singular} inspires the following definition, useful for the proof of Theorem~\ref{THM:near-dichotomy-for-cubic-curves-and-threefolds}.

\begin{definition}
\label{DEFN:special-singularity-of-cubic-threefold}
Let $X\belongs \PP^4_K$ be a $(3)$-CI, where $K=\ol{K}$.
Call $v\in X_{\map{sing}}(K)$ \emph{special} if $X$ contains a cone $S\belongs \PP^4_K$ with vertex $v$, where $S$ is a plane or a singular cubic $2$-scroll.
\end{definition}

\begin{proposition}
\label{PROP:uniform-characterization-of-low-degree-CI's-and-cubic-scrolls}
Let $Z$ be an integral closed subscheme of $\PP^n$ over an algebraically closed field $K$, where $n\geq 2$.
Then the following are equivalent:
\begin{enumerate}
    \item $\dim{Z} = n-2$, and $Z$ contains a nonempty open subscheme of a $(2,2)$-CI in $\PP^n$;
    \item $\abs{Z}$ is an irreducible component of a $(2,2)$-CI in $\PP^n$;
    and
    \item $Z$ is either a $(d_1,d_2)$-CI with $1\leq d_1\leq d_2\leq 2$, or a cubic $(n-2)$-scroll.
\end{enumerate}
\end{proposition}

\begin{proof}
The implication (1)$\Rightarrow$(2) is clear, since $Z$ is irreducible.

(2)$\Rightarrow$(3):
Assume (2); let the $(2,2)$-CI be $Y$.
Then $\deg{Z}\leq \deg{Y} = 4$.
If $\deg{Z}\leq 3$, then (3) follows from Lemma~\ref{LEM:description-of-low-degree-components-of-a-CI} (since $Z$ is integral).
If $\deg{Z}=4$, then a calculation using Hilbert polynomials shows that the inclusion morphism $Z\inject Y$ is birational, and $\abs{Z}=\abs{Y}$ (since $Y$ is equidimensional);
so $Y$ is reduced by Proposition~\ref{PROP:locally-Noetherian-generically-reduced-and-CM-implies-reduced} (since $Z$ is reduced), and therefore $Z=Y$;
so $Z$ is a $(2,2)$-CI, and (3) again holds.

(3)$\Rightarrow$(1):
If $Z$ is a cubic $(n-2)$-scroll, then Proposition~\ref{PROP:basic-cubic-scroll-facts}(\ref{ITEM:cubic-scroll-contains-a-nonempty-open-subscheme-of-a-(2,2)-CI}) implies (1).
Now suppose $Z$ is a $\bm{d}$-CI of the form $V(P_1,P_2)$, with $\bm{d}=(d_1,d_2)\in \set{1, 2} \times \set{1, 2}$.
Then $(P_1),(P_2)$ are coprime.
So if $L_1,L_2$ are sufficiently general linear forms, then $Y\defeq V(L_1^{2-d_1}P_1, L_2^{2-d_2}P_2)$ is a $(2,2)$-CI that coincides with $Z$ away from $V(L_1L_2)$.
So (1) again holds.
\end{proof}

To treat some lacunary cases, we need a final technical ingredient (Lemma~\ref{LEM:lower-bound-on-E-of-stably-non-conical-cubic-surfaces}).

\begin{lemma}
\label{LEM:lower-bound-on-E-of-stably-non-conical-cubic-surfaces}
Let $X\belongs \PP^3$ be a $(3)$-CI over a finite field $k$.
Suppose $X_{\ol{k}}$ is not a cone.
Then there exists a finite extension $k'/k$ such that $E(X_{k''})$ is either constant, or asymptotic to a nonconstant polynomial in $\card{k''}$ with positive leading coefficient, as we let $k''$ range over finite extensions of $k'$.
The latter holds if $X$ is smooth.
\end{lemma}

\begin{proof}
If $X$ is smooth, then it is known that for a suitable choice of $k'/k$, we have $E(X_{k''}) = 6\card{k''}$ (for all $k''/k'$).
Now suppose $X$ is singular.
By base change and a linear change of coordinates, assume $[0:0:0:1]\in X_{\map{sing}}$, and let $f_2,f_3$ be as in Proposition~\ref{PROP:using-k-rational-singular-point-to-count-points}.
Since $X_{\ol{k}}$ is not a cone, $f_2\neq 0$.
If $f_2,f_3$ share a nonconstant factor, then $X_{\ol{k}} = V_{\PP^3}(f_2x_4+f_3)_{/\ol{k}}$ is non-integral, and thus reducible (since $X_{\ol{k}}$ is not a cone), so (for a suitable choice of $k'/k$) Lang--Weil gives $E(X_{k''})\sim c\card{k''}^2$ for some $c>0$.

Now suppose $f_2,f_3$ are nonzero and coprime.
Eq.~\eqref{EQN:rationality-based-equation-for-E(X)} may be rearranged to get
\begin{equation*}
E(X_{k''}) = 1 + \card{k''}\cdot \card{V_{\PP^2}(f_2,f_3)(k'')} - \card{V_{\PP^2}(f_2)(k'')}.
\end{equation*}
But for a suitable choice of $k'/k$, both of the following hold:
\begin{itemize}
    \item $V_{\PP^2}(f_2,f_3)(k'')$ is a finite nonempty set independent of $k''$;
    and
    \item $\card{V_{\PP^2}(f_2)(k'')} = c\card{k''}+1$, where $c\in \set{1,2}$ denotes the number of distinct irreducible components of $V_{\PP^2}(f_2)_{/\ol{k}}$.
\end{itemize}
So we are done if we can rule out the case in which $\card{V_{\PP^2}(f_2,f_3)(\ol{k})}=1$ and $c=2$.
Suppose for contradiction that this is the case; since $c=2$, we may assume $f_2 = x_2x_3$ (after base change).
Then $\card{V_{\PP^2}(f_2,f_3)(\ol{k})}=1$ forces
\begin{equation*}
V_{\PP^2}(x_i,f_3)(\ol{k})\belongs V_{\PP^2}(x_2,x_3)(\ol{k}) = \set{[1:0:0]}
\end{equation*}
for all $i\in \set{2,3}$.
Thus $f_3\in (x_2,x_3^3)\cap (x_2^3,x_3) = (x_2x_3,x_2^3,x_3^3)$, whence $X_{\ol{k}} = V_{\PP^3}(x_2x_3x_4+f_3)_{/\ol{k}}$ is a cone, contradicting our assumptions.
\end{proof}

\begin{proof}
[Proof of Theorem~\ref{THM:near-dichotomy-for-cubic-curves-and-threefolds}]
First say $n=2$.
Then by Lang--Weil, $X$ is $\abs{E}$-good if and only if $X_{\ol{k}}$ is irreducible.
But $X_{\ol{k}}$ is not a cone over a cone, so $X_{\ol{k}}$ is reducible if and only if it contains a line in $\PP^2_{\ol{k}}$.
So $X$ is $\abs{E}$-bad if and only if $X_{\ol{k}}$ contains a line in $\PP^2_{\ol{k}}$.

Now say $n=4$.
We will first show that conditions~(\ref{ITEM:cubic-threefold-stably-E-bad})--(\ref{ITEM:cubic-threefold-stably-|E|-bad}) are equivalent, and that they imply condition~(\ref{ITEM:cubic-threefold-potentially-contains-a-special-subvariety}); and that conditions~(\ref{ITEM:cubic-threefold-stably-E-bad})--(\ref{ITEM:cubic-threefold-potentially-contains-a-special-subvariety}) are equivalent if $\dim(X_{\map{sing}})\leq 0$.
If $X$ is smooth, then conditions~(\ref{ITEM:cubic-threefold-stably-E-bad})--(\ref{ITEM:cubic-threefold-stably-|E|-bad}) fail by the Weil conjectures (for our $X$, see \cites{bombieri1967local,manin1968correspondences} for proofs in the spirit of the present work),
and condition~(\ref{ITEM:cubic-threefold-potentially-contains-a-special-subvariety}) fails by Lemma~\ref{LEM:any-special-subvariety-is-a-cone-with-vertex-singular}.
Now assume $X$ is singular.
By base change (via Lemma~\ref{LEM:potential-and-stable-logic}) and a linear change of coordinates,
we may assume that $[0:0:0:0:1]\in X_{\map{sing}}$,
and that $[0:0:0:0:1]$ is a \emph{special} singularity of $X_{\ol{k}}$ (see Definition~\ref{DEFN:special-singularity-of-cubic-threefold}) if such a singularity exists.
Let $f_2,f_3$ be as in Proposition~\ref{PROP:using-k-rational-singular-point-to-count-points}.

\emph{Case~1: $f_2=0$.}
Then $f_3\neq 0$,
whence $X$ is a cone over a (possibly singular) cubic hypersurface $W\belongs V_{\PP^4}(x_5)\cong \PP^3$,
and therefore $X_{\ol{k}}$ contains a plane in $\PP^4_{\ol{k}}$ (since it is known that $W_{\ol{k}}$ must contain a line; see e.g.~\cite{mustata2017lecture2}*{Theorem~1.1}).
In particular, condition~(\ref{ITEM:cubic-threefold-potentially-contains-a-special-subvariety}) holds.
Also, by Proposition~\ref{PROP:cone-construction-multiplies-E-by-q} and Lemma~\ref{LEM:lower-bound-on-E-of-stably-non-conical-cubic-surfaces}, conditions~(\ref{ITEM:cubic-threefold-stably-E-bad})--(\ref{ITEM:cubic-threefold-stably-|E|-bad}) are equivalent; and if $\dim(X_{\map{sing}})\leq 0$, then $W$ is smooth, so conditions~(\ref{ITEM:cubic-threefold-stably-E-bad})--(\ref{ITEM:cubic-threefold-stably-|E|-bad}) hold.

\emph{Case~2: $\dim V_{\PP^3}(f_2,f_3) = \dim V_{\PP^3}(f_2)$.}
Then $f_2\neq 0$, and $f_2,f_3$ share a nonconstant factor.
So the cubic $f_2x_5+f_3$ is reducible, whence $X_{\ol{k}}$ contains a $3$-plane in $\PP^4_{\ol{k}}$.
But $X_{\ol{k}}$ is not a cone over a cone, so $X_{\ol{k}}$ is reducible, and thus conditions~(\ref{ITEM:cubic-threefold-stably-E-bad})--(\ref{ITEM:cubic-threefold-potentially-contains-a-special-subvariety}) all hold.

\emph{Case~3: $f_2\neq 0$ and $\dim V_{\PP^3}(f_2,f_3)\neq \dim V_{\PP^3}(f_2)$.}
In other words, $f_2,f_3$ are nonzero and coprime.
So by Lemma~\ref{LEM:associated-(2,3)-CI-controls-potential-goodness-of-singular-cubic-threefold},
conditions~(\ref{ITEM:cubic-threefold-stably-E-bad})--(\ref{ITEM:cubic-threefold-stably-|E|-bad}) are equivalent,
and they imply that $V_{\PP^3}(f_2,f_3)_{/\ol{k}}$ contains a $(1,1)$-CI, a $(1,2)$-CI, or a cubic $1$-scroll, in $\PP^3_{\ol{k}}$;
and this in turn implies condition~(\ref{ITEM:cubic-threefold-potentially-contains-a-special-subvariety}), by Lemma~\ref{LEM:convert-(1,2)-CI-or-(2,2)-CI-containment-to-(1,1)-CI-containment} and Proposition~\ref{PROP:basic-cubic-scroll-facts}(\ref{ITEM:singular-cubic-scroll-equivalences}) (since $V_{\PP^4}(f_2,f_3)\belongs X$).
Conversely, suppose condition~(\ref{ITEM:cubic-threefold-potentially-contains-a-special-subvariety}) holds.
Then by Lemma~\ref{LEM:any-special-subvariety-is-a-cone-with-vertex-singular}, $X_{\ol{k}}$ has a special singularity.
Thus (by assumption) $[0:0:0:0:1]$ is special,
so (by Definition~\ref{DEFN:special-singularity-of-cubic-threefold} and the basic properties of cones) $X_{\ol{k}}$ contains the cone $C_{[0:0:0:0:1]}(Y)\belongs \PP^4_{\ol{k}}$ over some $Y\belongs V_{\PP^4}(x_5)_{/\ol{k}}$, where $Y$ is either a $(1,1)$-CI or a cubic $1$-scroll.
The proof of Lemma~\ref{LEM:any-special-subvariety-is-a-cone-with-vertex-singular} then gives $Y\belongs V_{\PP^3}(f_2,f_3)_{/\ol{k}}$.
So if $\dim(X_{\map{sing}})\leq 0$, then conditions~(\ref{ITEM:cubic-threefold-stably-E-bad})--(\ref{ITEM:cubic-threefold-stably-|E|-bad}) hold by Lemma~\ref{LEM:associated-(2,3)-CI-controls-potential-goodness-of-singular-cubic-threefold}, as desired.

To complete the proof of Theorem~\ref{THM:near-dichotomy-for-cubic-curves-and-threefolds}, first note that a $(2,2)$-CI of the form $V(Q_1,Q_2)\belongs \PP^4_{\ol{k}}$ has $V(Q_1)_{\map{sing}} \cap V(Q_2)_{\map{sing}} \neq \emptyset$ if and only if $V(Q_1,Q_2)$ is a cone.

(\ref{ITEM:cubic-threefold-potentially-contains-a-special-subvariety})$\Rightarrow$(\ref{ITEM:cubic-threefold-potentially-contains-a-nonempty-open-subscheme-of-a-special-CI}):
Assume condition~(\ref{ITEM:cubic-threefold-potentially-contains-a-special-subvariety}) holds.
Then by Proposition~\ref{PROP:basic-cubic-scroll-facts}(\ref{ITEM:singular-cubic-scroll-equivalences}), $X_{\ol{k}}$ contains a cone of the form $C_v(Z)\belongs \PP^4_{\ol{k}}$, where $Z$ is either a line, or a cubic $1$-scroll, in some $3$-plane $H\belongs \PP^4_{\ol{k}}$ with $v\notin H$.
Now apply Proposition~\ref{PROP:uniform-characterization-of-low-degree-CI's-and-cubic-scrolls}(3)$\Rightarrow$(1) to $Z\belongs H$.

(\ref{ITEM:cubic-threefold-potentially-contains-a-nonempty-open-subscheme-of-a-special-CI})$\Rightarrow$(\ref{ITEM:cubic-threefold-potentially-contains-a-special-subvariety}):
Assume condition~(\ref{ITEM:cubic-threefold-potentially-contains-a-nonempty-open-subscheme-of-a-special-CI}) holds.
Then by Proposition~\ref{PROP:uniform-characterization-of-low-degree-CI's-and-cubic-scrolls}(2)$\Rightarrow$(3), $X_{\ol{k}}$ contains a cone in $\PP^4_{\ol{k}}$ that is either (i) a $\bm{d}$-CI for some $\bm{d}\in \set{1, 2} \times \set{1, 2}$; or (ii) a cone over a cubic $1$-scroll.
In case~(i), use Lemma~\ref{LEM:convert-(1,2)-CI-or-(2,2)-CI-containment-to-(1,1)-CI-containment}; in case~(ii), use Proposition~\ref{PROP:basic-cubic-scroll-facts}(\ref{ITEM:singular-cubic-scroll-equivalences}).
\end{proof}

\section{Analysis of the Fermat cubic fourfold}
\label{SEC:Fermat-cubic-fourfold-analysis}

The main goal of this section is to prove the rather specialized Theorem~\ref{THM:main-diagonal-theorem}; but let us first give some more general context.
Consider the following question motivated by
Theorem~\ref{THM:main-special-result-informal}, Problem~\ref{PROB:moments-of-failures-of-sqrt-cancellation}, and Theorem~\ref{THM:concrete-moment-and-super-level-bounds}(3):

\begin{question}
Given a smooth cubic hypersurface $X\belongs \PP^5_\CC$, let $S\belongs (\PP^5_\CC)^\vee(\CC)$ denote (via the interpretation of $(\PP^5_\CC)^\vee$ as a Grassmannian) the set of hyperplanes $H\belongs \PP^5_\CC$ such that $X\cap H$ contains a singular cubic $2$-scroll in $\PP^5_\CC$.
What is the best possible upper bound on the dimension of the closure $\ol{S}\belongs (\PP^5_\CC)^\vee$?
\end{question}

\begin{remark}
It is known that $\#S = 1$ (and thus $\dim\ol{S} = 0$) for sufficiently general $X$ containing a (not necessarily singular) cubic $2$-scroll in $\PP^5_\CC$; cf.~especially \cite{hassett1996special}*{proof of Lemma~2.11, and the subsequent dimension counting} and \cite{hassett2010flops}*{Propositions~3.3 and~6.1}.

Also, $S$ may well be finite for $X = V_{\PP^5}(x_1^3+\dots+x_6^3)_{/\CC}$;
our analysis below in \S\ref{SEC:Fermat-cubic-fourfold-analysis} falls short of a proof (due to interference from the planes on $X$),
but Proposition~\ref{PROP:main-result-on-singular-cubic-scrolls-in-Fermat-cubic-fourfold} (the key to Theorem~\ref{THM:main-diagonal-theorem}) does imply for $X = V_{\PP^5}(x_1^3+\dots+x_6^3)_{/\CC}$ that $\dim\ol{S}$ is at most $2$ (the dimension of the space of hyperplanes $H$ such that $X\cap H$ contains a plane in $\PP^5_\CC$).

In any case, planes ``generically'' seem to play a larger role than singular cubic $2$-scrolls, so it seems quite reasonable to believe $\dim\ol{S}\leq 2$ (if not a stronger bound).
If $\dim\ol{S}\leq 2$ is indeed true in general, and if a suitable analog remains true in positive characteristic, then one would likely be able to remove the diagonality assumption on $G$ in Theorem~\ref{THM:concrete-moment-and-super-level-bounds}(3).
\end{remark}


\begin{proof}
[Proof of Theorem~\ref{THM:main-diagonal-theorem}, assuming Proposition~\ref{PROP:main-result-on-singular-cubic-scrolls-in-Fermat-cubic-fourfold} below]
Let $m,F_1,\dots,F_m,k,F$ be as in Theorem~\ref{THM:main-diagonal-theorem}.
Assume the characteristic of $k$ is sufficiently large in terms of the integers $F_1,\dots,F_m$.
If $V_{\bm{c}}$ fails ``square-root cancellation'', then by Theorem~\ref{THM:main-special-result-informal} $(V_{\bm{c}})_{\ol{k}}$ contains either an $\frac{m-2}{2}$-plane or a singular cubic $2$-scroll, and by Proposition~\ref{PROP:main-result-on-singular-cubic-scrolls-in-Fermat-cubic-fourfold} we may assume that $(V_{\bm{c}})_{\ol{k}}$ contains an $\frac{m-2}{2}$-plane.
Furthermore, the $\frac{m-2}{2}$-planes on $V_{\ol{k}}$ are known to be cut out by systems of equations of the form ``$F_ix_i^3+F_jx_j^3=0$ in pairs'' (see e.g.~\cite{wang2022thesis}*{Remark~6.3.8 and references within}).
Thus $V_{\bm{c}}$ fails ``square-root cancellation'' if and only if ``$c_i^3/F_i=c_j^3/F_j$ in pairs''.
\end{proof}

For the rest of \S\ref{SEC:Fermat-cubic-fourfold-analysis}, let $K$ be an algebraically closed field with characteristic $p\geq 0$, and let $F\defeq x_1^3+\dots+x_6^3$.
Define $V,V_{\bm{c}}$ following the Main Setup (taking $k=K$ and $m=6$).
For $\bm{c}=(c_1,\dots,c_6)\in K^6\setminus \set{\bm{0}}$ we let $[\bm{c}]^\perp=[c_1:\dots:c_6]^\perp$ denote the $5$-dimensional subspace
\begin{equation*}
\set{\bm{x}=(x_1,\dots,x_6)\in K^6: c_1x_1+\dots+c_6x_6=0}
\end{equation*}
of $K^6$,
and for a vector space $M$ over $K$ we let $\PP(M)$ denote the \emph{projectivization} of $M$.

Let $[\bm{a}]=[a_1:\dots:a_6]\in V(K)$.
Let $\ell_{[\bm{a}]}$ denote the one-dimensional subspace ``$K\cdot \bm{a}$'' of $K^6$.
To analyze singular scrolls on $V$, we now define the \emph{tangential vision} $Y=Y^{[\bm{a}]}$ of $[\bm{a}]$.

\begin{definition}
Let $Y$ denote the closed subscheme of $\PP(K^6/\ell_{[\bm{a}]})\cong \PP^4$ defined, upon choosing a $6\times 5$ matrix $B$ whose column vectors span $K^6/\ell_{[\bm{a}]}$, by the homogeneous ideal of $K[\bm{y}]=K[y_1,\dots,y_5]$ generated by the coefficients of the cubic $F(s\bm{a}+tB\bm{y})\in (K[\bm{y}])[s,t]$.
\end{definition}

\begin{remark}
\label{RMK:tangential-vision-calculations}
Canonically, $Y\belongs \PP([\grad{F}(\bm{a})]^\perp/\ell_{[\bm{a}]})$ and $C_{[\bm{a}]}(Y) = V_{\PP^5}(F(s\bm{a}+t\bm{x}))\belongs V_{[\grad{F}(\bm{a})]}$.
Now suppose $a_6\neq 0$, so that $K^5\times \set{0}$ spans $K^6/\ell_{[\bm{a}]}$.
Write $F(s\bm{a}+(t\bm{y},0)) = f_1s^2t+f_2st^2+f_3t^3\in (K[\bm{y}])[s,t]$; then we may identify $Y$ with $V_{\PP^4}(f_1,f_2,f_3)$.
A short calculation then shows that $V_{[\grad{F}(\bm{a})]}\cong V_{\PP^5}(f_2y_6+f_3, f_1)$ (if we identify $\bm{x}$ with $y_6\bm{a}+(\bm{y},0)$).
\end{remark}

\begin{proposition}
\label{PROP:main-result-on-singular-cubic-scrolls-in-Fermat-cubic-fourfold}
Let $[\bm{a}]\in V(K)$.
The following conditions are equivalent:
\begin{enumerate}
    \item $V$ contains a singular cubic $2$-scroll $S\belongs \PP^5$ with $[\bm{a}]\in S_{\map{sing}}$;
    and
    \item $Y=Y^{[\bm{a}]}$ contains a cubic $1$-scroll in $\PP(K^6/\ell_{[\bm{a}]})$.
\end{enumerate}
Now \emph{suppose} (1) holds for some $S$.
Then $\Span{S} = \PP([\grad{F}(\bm{a})]^\perp) \belongs \PP^5$.
Furthermore, there exists a constant $C>0$ such that if $V\cap (\Span{S})$ contains no plane in $\PP^5$, then $0 < p\le C$.
\end{proposition}

The proof needs two lemmas.
The basic idea is that a union of two twisted cubic curves is singular at each intersection point, and that the degree of the intersection can be bounded from below using Hilbert polynomials.
On the other hand, we can bound the singularities from above by computer calculation, after first controlling their severity by a combinatorial singularity analysis based on diagonality and Vandermonde determinants.
With more work, one might be able to determine precisely when the condition \ref{PROP:main-result-on-singular-cubic-scrolls-in-Fermat-cubic-fourfold}(1) holds.
The main bottlenecks in the present argument are Lemmas~\ref{LEM:conditions-on-vision-imposed-by-twisted-cubic}(4) and~\ref{LEM:linearly-rich-configurations}.

\begin{lemma}
\label{LEM:conditions-on-vision-imposed-by-twisted-cubic}
Suppose $p\nmid 6$.
Let $[\bm{a}]\in V(K)$.
Let $n\defeq -1 + \#\set{i\in [6]: a_i\neq 0}$.
Assume $Y=Y^{[\bm{a}]}$ contains a cubic $1$-scroll in $\PP(K^6/\ell_{[\bm{a}]})$.
If $n\neq 1$, then the following hold:
\begin{enumerate}
    \item $n\in \set{4,5}$ and $\dim{Y} = 1$, and $\dim(Y_{\map{sing}}) \leq 0$.
    \item Locally near each of its points, $Y_{\map{sing}}$ has degree exactly $6-n$ in $\PP(K^6/\ell_{[\bm{a}]})$.
    \item $\#Y_{\map{sing}}(K) = -1 + 2^{n-6}\cdot \#\set{\textnormal{subsets}\;I\belongs [6]: \sum_{i\in I} a_i^3 = 0}$.
    \item If $Y$ has no irreducible component of degree $\leq 2$, then $Y_{\map{sing}}$ has degree $\geq 5$.
    \item If $[a_1^3:\dots:a_6^3] = [-2: -1: 1: 1: 0: 1]$,
    then $0<p\le \mcal{B}$ (for some constant $\mcal{B}$).
\end{enumerate}
\end{lemma}

Below, let $\map{LinAut}(V)$ denote the group of (linear) automorphisms of $\PP^5$ that preserve $V$.

\begin{proof}
Via $\map{LinAut}(V)$, we may assume $a_6\neq 0$.
Now let $f_1,f_2,f_3$ be as in Remark~\ref{RMK:tangential-vision-calculations}; explicitly, $(f_1,f_2,f_3) = (3a_1^2y_1+\dots+3a_5^2y_5, 3a_1y_1^2+\dots+3a_5y_5^2, y_1^3+\dots+y_5^3)$.

Note that $n\geq 1$, since $F(\bm{a})=0$ and $\bm{a}\neq \bm{0}$.
So we may assume $n\geq 2$.
A short calculation then shows that $f_1\nmid f_2$ (the only subtle case is $n=2$; if $a_3=a_4=a_5=0$ and $a_1a_2\neq 0$, say, then we need to use the fact that $a_1^3+a_2^3 = -a_6^3\neq 0$).
Now write $K[\bm{y}]/(f_1)\cong K[z_1,z_2,z_3,z_4]$, and let $g_2,g_3\in K[\bm{z}]$ correspond to $f_2,f_3\bmod{f_1}$; then $g_2\neq 0$.
Furthermore, $V_{[\grad{F}(\bm{a})]}$ is irreducible by Theorem~\ref{THM:Zak's-principle}, so $g_2,g_3$ are coprime.
Thus $\dim{Y} = 1$, i.e.~$Y$ is a $(2,3)$-CI in $V_{\PP^4}(f_1)\cong \PP^3$.
By assumption, $Y$ contains a cubic $1$-scroll.
By Lemma~\ref{LEM:description-of-low-degree-components-of-a-CI}, it follows that $g_2$ is irreducible.
In particular, we must have $n-1\geq 3$, so $n\in \set{4,5}$.
Furthermore, $\dim(V_{\PP^3}(g_2)_{\map{sing}})\leq 0$, so by Theorem~\ref{THM:Zak's-principle} and Lemma~\ref{LEM:lower-bound-on-singular-locus-in-terms-of-any-associated-(2,3)-CI}(1), $\dim(Y_{\map{sing}}) \leq 0$.
This completes the proof of (1).

Now we study $Y_{\map{sing}}$.
By the definitions, $Y_{\map{sing}}$ is defined in $Y$ by the $3\times 3$ minors of the $3\times 5$ matrix with rows $(3a_1^2,\dots,3a_5^2)$, $(6a_1y_1,\dots,6a_5y_5)$, $(3y_1^2,\dots,3y_5^2)$.
Vandermonde thus gives $Y_{\map{sing}} = V_{\PP^4}(f_1,f_2,f_3,(a_iy_j-a_jy_i)(a_iy_k-a_ky_i)(a_jy_k-a_ky_j))_{\set{i,j,k}}$, where $\set{i,j,k}$ ranges over $\binom{[5]}{3}$ (the set of unordered triples in $[5]$).
Let $U\defeq \set{i\in [5]: a_i\neq 0}$ (so $\card{U} = n$).

Let $[\bm{z}]\in Y_{\map{sing}}(K)$; then $\#\set{z_i/a_i: i\in U}\leq 2$.
But also $f_1=0$ and $a_1^3+\dots+a_5^3 = -a_6^3\neq 0$, so $\#\set{z_i/a_i: i\in U}\geq 2$.
Thus there exist distinct $r_1,r_2\in K$ such that the sets $I_l\defeq \set{i\in U: z_i/a_i = r_l}$ for $l\in [2]$ are nonempty, and partition $U$.
A short calculation then shows that locally near $[\bm{z}]$, the scheme $Y_{\map{sing}}$ coincides with
\begin{equation*}
\begin{split}
&V_{\PP^4}(f_1,f_2,f_3,
(a_iy_j-a_jy_i)_{\set{i,j}\belongs I_1},
(a_iy_j-a_jy_i)_{\set{i,j}\belongs I_2},
(y_k^2)_{k\in [5]\setminus U}) \\
&\cong V_{\PP^{6-n}}(A_1w_1+A_2w_2, A_1w_1^2+A_2w_2^2, A_1w_1^3+A_2w_2^3, w_3^2\cdot \bm{1}_{n=4})\eqdef W,
\end{split}
\end{equation*}
where $A_l\defeq \sum_{i\in I_l} a_i^3$, and the variable $w_l$ (for $l\in [2]$) corresponds to $y_i/a_i$ for each $i\in I_l$.
Here $A_1+A_2 = -a_6^3\neq 0$, and $[\bm{z}]$ corresponds to a point in $W$ (whence $W\neq \emptyset$), so a short calculation shows that $A_1A_2=0$, and if $A_l=0$, then $A_{3-l}\neq 0$ and $r_{3-l}=0\neq r_l$.
Furthermore, $\#W(K) = 1$ and $\deg{W} = 6-n$.
This completes the proof of (2).

Conversely, if $\sum_{i\in I} a_i^3 = 0$ for some nonempty set $I\belongs U$, then a short calculation shows that $V_{\PP^4}((a_iy_j-a_jy_i)_{\set{i,j}\belongs I}, (y_k)_{k\in [5]\setminus I})\belongs Y_{\map{sing}}$.
This implies (3).

Since $g_2$ is irreducible, each degree-$3$ irreducible component of $Y$ is (by Lemma~\ref{LEM:description-of-low-degree-components-of-a-CI}) a cubic $1$-scroll in $V_{\PP^4}(f_1)\cong \PP^3$.
And $Y$ is reduced by Proposition~\ref{PROP:locally-Noetherian-generically-reduced-and-CM-implies-reduced}, since $\dim(Y_{\map{sing}}) \leq 0$.
For the proof of (4) only, suppose $Y$ has no irreducible component of degree $\leq 2$; then $Y$ must be a scheme-theoretic union of two cubic $1$-scrolls in $V_{\PP^4}(f_1)$, say $T_1,T_2$, with $\dim(T_1\cap T_2) = 0$.
But a $(2,3)$-CI in $V_{\PP^4}(f_1)$ has Hilbert polynomial $6d-3$, and a cubic $1$-scroll has Hilbert polynomial $3d+1$, so $\deg(T_1\cap T_2) = 2(3d+1) - (6d-3) = 5$.
Furthermore, $Y$ is contained in surfaces $Z_2\defeq V_{\PP^4}(f_1,f_2)$ and $Z_3\defeq V_{\PP^4}(f_1,f_3)$ with disjoint singular loci (a short calculation, using $a_1^3+\dots+a_5^3\neq 0$, shows that $Z_2$ is smooth if $n=5$, and singular at the point $[\bm{1}_{i\notin U}]_{i\in [5]}$ if $n=4$; but $[\bm{1}_{i\notin U}]_{i\in [5]}\notin Z_3$ if $n=4$),
so $Y_{\map{sing}}\contains T_1\cap T_2$ by Proposition~\ref{PROP:lower-bound-on-singular-scheme-of-a-reducible-variety}.
Thus $\deg(Y_{\map{sing}})\geq 5$, proving (4).

Finally, if $[a_1^3:\dots:a_6^3] = [-2: -1: 1: 1: 0: 1]$ and $p=0$,
then by the Magma calculation
\begin{verbatim}
R<t> := PolynomialRing(Rationals()); K<a> := NumberField(t^3-2);
P<y_1, y_2, y_3, y_4, y_5> := ProjectiveSpace(K, 4);
C := Curve(P, [a^2*y_1+y_2+y_3+y_4, -a*y_1^2-y_2^2+y_3^2+y_4^2,
y_1^3+y_2^3+y_3^3+y_4^3+y_5^3]); [IsReduced(C), IsAbsolutelyIrreducible(C)];
\end{verbatim}
(which outputs \verb![ true, true ]!)
and \cite{stacks-project}*{\href{https://stacks.math.columbia.edu/tag/020I}{Tag 020I}} (which implies that over $\QQ(2^{1/3})$, ``reduced'' implies ``geometrically reduced''),
$Y$ is geometrically integral,
since $Y\in \map{LinAut}(V)\cdot C$.
Since geometric integrality spreads out (see e.g.~\cite{poonen2017rational}*{Theorem~3.2.1}),
(5) follows.
\end{proof}

\begin{lemma}
\label{LEM:linearly-rich-configurations}
Suppose $p$ is $0$ or sufficiently large.
Let $\bm{x}\in K^6$ and $n\in \set{4,5}$.
Suppose $\sum_{i\in [6]} x_i = 0$ and $\#\set{i\in [6]: x_i=0} = 5-n$.
Assume $-1 + 2^{n-6}\cdot \#\set{I\belongs [6]: \sum_{i\in I} x_i = 0}\geq 3\cdot \bm{1}_{n=4} + 5\cdot \bm{1}_{n=5}$.
Then there exists $\pi\in S_6$ such that $[x_{\pi(1)}:\dots:x_{\pi(6)}]$ is $[-2: -1: 1: 1: 0: 1]$, $[-2: -2: 1: 1: 1: 1]$, or $[x_{\pi(1)}:-x_{\pi(1)}:-x_{\pi(1)}:x_{\pi(1)}:-1:1]$.
\end{lemma}

\begin{proof}
We may assume that $x_6=1$ and $x_1\cdots x_n\neq 0$.
Then by the hypotheses,
$3\cdot \bm{1}_{n=4} + 5\cdot \bm{1}_{n=5}
\leq -1 + 2^{n-6}\cdot \#\set{I\belongs [6]: \sum_{i\in I} x_i = 0}
= \#\set{\emptyset\neq I\belongs [n]: \sum_{i\in I} x_i = 0}$,
and $x_1+\dots+x_n = -1\neq 0$.
To classify all possible $\bm{x}$,
it remains to do a large but finite amount of linear algebra, which we do with computer assistance.
The SageMath code
\begin{verbatim}
def PermTest(N, S1, S2): # Think of S1, S2 as sets of sets of indices.
    for P in Permutations(range(N)):
        test = True
        for I in S1: PI = [P[i] for i in I]; PI.sort(); test = test and (PI in S2)
        if test and (len(S1) == len(S2)): return True
    return False
def RichSols(n,r,K):
    C = []
    for k in range(1,n+1): C.extend(Combinations(range(n),k).list())
    T = Combinations(C,r).list(); SimpleRichS = []; OtherRichS = []
    for S in T:
        A = matrix(K, r+1, n)
        for i in range(r):
            for j in S[i]: A[i,j] = 1
        rk = rank(A); bad = False
        for j in range(n): Aj = copy(A); Aj[r,j] = 1; bad = bad or (rank(Aj) == rk)
        for j in range(n): A[r,j] = 1
        if bad or (rank(A) == rk): continue
        if rank(A) == n:
            b = vector(K, r); b = vector(K, b.list() + [-1])
            v = list(A.solve_right(b)); v.sort()
            if v not in SimpleRichS: SimpleRichS.append(v)
        elif all([not(PermTest(n,S,U)) for U in OtherRichS]): OtherRichS.append(S)
    return SimpleRichS + ["---"] + OtherRichS
\end{verbatim}
consists of two functions.
{\tt PermTest} tests if two sets of sets $I\belongs [N]$ are identical up to $S_N$.
The main function {\tt RichSols} considers sets $S$ of $r$ nonempty sets $I\belongs [n]$;
discards those $S$ for which the scheme $L_S\defeq \bigcap_{I\in S} V_{\Aff^n}(\sum_{i\in I} x_i)_{/K}$ is contained in $V_{\Aff^n}(x_1\cdots x_n(x_1+\dots+x_n))_{/K}$;
and for the remaining $S$,
outputs (up to symmetry) each intersection $v_S\defeq L_S\cap (x_1+\dots+x_n=-1)$ of size one,
and outputs (up to symmetry) each $S$ for which $\dim_K(v_S) \geq 1$.

We now establish the conclusion of Lemma~\ref{LEM:linearly-rich-configurations}.
First suppose $p=0$.
The output
\begin{verbatim}
[[-1/2, -1/2, -1/2, 1/2], [-2, -1, 1, 1], [-1, -1, -1, 2], '---']
\end{verbatim}
of {\tt RichSols(4,3,QQ)} and the output
\begin{verbatim}
[[-1, -1, -1, 1, 1], [-1, -1/2, -1/2, 1/2, 1/2], [-2, -1, -1, 1, 2],
 [-1/2, -1/2, -1/2, -1/2, 1], [-2, -2, 1, 1, 1], '---',
 [[0, 1], [0, 2], [1, 3], [2, 3], [0, 1, 2, 3]],
 [[0, 1], [0, 2], [3, 4], [0, 1, 3, 4], [0, 2, 3, 4]]]
\end{verbatim}
of {\tt RichSols(5,5,QQ)} imply the desired results for $n=4$ and $n=5$, respectively.
For $n=4$, note that $[-\frac12: -\frac12: -\frac12: \frac12: 0: 1]$ and $[-1: -1: -1: 2: 0: 1]$
fall under the case $[-2: -1: 1: 1: 0: 1]$.
For $n=5$, note that $[-1: -1: -1: 1: 1: 1]$ and $[-1: -\frac12: -\frac12: \frac12: \frac12: 1]$
fall under the case $[x_{\pi(1)}:-x_{\pi(1)}:-x_{\pi(1)}:x_{\pi(1)}:-1:1]$,
as do $[-2: -1: -1: 1: 2: 1]$ (after scaling) and $[1:-1:-1:x_4:-x_4:1]$ (since $x_4\neq 0$);
and that $[-\frac12: -\frac12: -\frac12: -\frac12: 1: 1]$
falls under the case $[-2: -2: 1: 1: 1: 1]$.

Now suppose $p>0$.
There exists a finite set of integer matrices $A$ such that the outputs of {\tt RichSols(4,3,K)} and {\tt RichSols(5,5,K)} depend only on the function $A\mapsto \rank_{\FF_p}(A)$.
But for a fixed $A$, the determinantal characterization of rank implies that $\rank_{\FF_p}(A) = \rank_{\QQ}(A)$ for all sufficiently large $p$.
So the proof for $p=0$ carries over (for sufficiently large $p$).
\end{proof}

\begin{proof}
[Proof of Proposition~\ref{PROP:main-result-on-singular-cubic-scrolls-in-Fermat-cubic-fourfold}]
By Remark~\ref{RMK:tangential-vision-calculations}, $C_{[\bm{a}]}(Y)\belongs V$.
A cone over a cubic $1$-scroll is a cubic $2$-scroll singular at the vertex of the cone (and smooth elsewhere); cf.~Proposition~\ref{PROP:basic-cubic-scroll-facts}(\ref{ITEM:singular-cubic-scroll-equivalences}).
So (2) implies (1).

Now suppose (1) holds for some $S$.
By Lemma~\ref{LEM:any-special-subvariety-is-a-cone-with-vertex-singular}, there exists $v\in (V\cap (\Span{S}))_{\map{sing}}(K)$ such that $S$ is a cone with vertex $v$.
Here we must have $\set{v} = S_{\map{sing}} = \set{[\bm{a}]}$ (by the structure theory of cubic scrolls).
So $S$ is a cone with vertex $[\bm{a}]$, whence $S\belongs C_{[\bm{a}]}(Y)$ by Remark~\ref{RMK:tangential-vision-calculations} and the proof of Lemma~\ref{LEM:any-special-subvariety-is-a-cone-with-vertex-singular}.
Thus (2) holds.

Furthermore, $V\cap (\Span{S})$ is singular at $v = [\bm{a}]$, so $\Span{S} = \PP([\grad{F}(\bm{a})]^\perp)$.
In particular, $V\cap (\Span{S}) = V_{[\grad{F}(\bm{a})]}$.
But by Remark~\ref{RMK:tangential-vision-calculations}, $V_{[\grad{F}(\bm{a})]}\contains C_{[\bm{a}]}(Y)$ and $Y\belongs \PP([\grad{F}(\bm{a})]^\perp/\ell_{[\bm{a}]})$.
Suppose $V_{[\grad{F}(\bm{a})]}$ contains no plane; then by Lemma~\ref{LEM:convert-(1,2)-CI-or-(2,2)-CI-containment-to-(1,1)-CI-containment}, it contains no $(1,1)$-CI or $(1,2)$-CI in $\PP([\grad{F}(\bm{a})]^\perp)$.
So $Y$ contains no $(1,1)$-CI or $(1,2)$-CI in $\PP([\grad{F}(\bm{a})]^\perp/\ell_{[\bm{a}]})$.
By Lemma~\ref{LEM:description-of-low-degree-components-of-a-CI}, it follows that either $\dim{Y}\geq 2$, or $Y$ is a $(2,3)$-CI curve in $\PP([\grad{F}(\bm{a})]^\perp/\ell_{[\bm{a}]})$ with no irreducible component of degree $\leq 2$.
We now proceed by casework on $n\defeq -1 + \#\set{i\in [6]: a_i\neq 0}$.

\emph{Case~1: $n=1$.}
Then $\sigma([\bm{a}]) = [0:0:0:0:-1:1]$ for some $\sigma\in \map{LinAut}(V)$.

\emph{Case~2: $n\neq 1$.}
If $p$ is either $0$, or sufficiently large (say $\geq A$) in terms of Lemmas~\ref{LEM:conditions-on-vision-imposed-by-twisted-cubic}(5) and~\ref{LEM:linearly-rich-configurations},
then by Lemmas~\ref{LEM:conditions-on-vision-imposed-by-twisted-cubic} and~\ref{LEM:linearly-rich-configurations}, we have $n=5$, and there exists $\pi\in S_6$ such that $[a_{\pi(1)}^3:\dots:a_{\pi(6)}^3]$ is $[-2: -2: 1: 1: 1: 1]$ or $[a_{\pi(1)}^3:-a_{\pi(1)}^3:-a_{\pi(1)}^3:a_{\pi(1)}^3:-1:1]$.

If $0<p<A$, we are done, so suppose otherwise.
Then whether $n=1$ or $n\neq 1$, there exist $\sigma\in \map{LinAut}(V)$ and $c\in K$ with $\sigma([\grad{F}(\bm{a})]^\perp) = [c:c:c:c:1:1]^\perp$.
So $V_{[\grad{F}(\bm{a})]}$ contains the plane $\sigma^{-1}(V(x_1+x_2, x_3+x_4, x_5+x_6))$.
This contradicts our assumption on $V_{[\grad{F}(\bm{a})]}$.
\end{proof}

\appendix

\section{Moment calculations}
\label{SEC:moment-calculations}

In this appendix, we prove two results about moments:
Theorem~\ref{THM:concrete-moment-and-super-level-bounds} (proven using Theorems~\ref{THM:main-general-result-informal} and~\ref{THM:main-diagonal-theorem})
and Corollary~\ref{COR:smooth-locus-moment-calculations} (to Theorem~\ref{THM:main-general-result-informal}).

\begin{proof}
[Proof of Theorem~\ref{THM:concrete-moment-and-super-level-bounds}]
Suppose $p$ is sufficiently large in terms of $G,d,m$.
Then $\abs{E_{\bm{c}}(q)}\leq C_1(d,m) q^{(m-2)/2}$ holds for all $\bm{c}\in \FF_q^m\setminus \set{\bm{0}}$;
this can be proven, for instance, using \cite{hooley1991number}*{Theorem~2} and the present Theorem~\ref{THM:Zak's-principle}.
Furthermore, by the Weil conjectures, $\abs{E_{\bm{c}}(q)}\leq C_2(d,m)q^{(m-3)/2}$ holds for all $\bm{c}\in \FF_q^m$ with $\disc(G,\bm{c})\neq 0$.
By Lang--Weil for $\disc(G,\bm{c})=0$, claim~(1) readily follows.
If $d\geq 3$, use Theorem~\ref{THM:main-general-result-informal} (plus Proposition~\ref{PROP:amplification-for-projective-varieties}) instead of the Weil conjectures; a similar argument then yields claim~(2).
If $(m,d) = (6,3)$ and $G$ is diagonal, similarly use Theorem~\ref{THM:main-diagonal-theorem} (plus Proposition~\ref{PROP:amplification-for-projective-varieties}) to prove claim~(3).
\end{proof}

It remains to prove Corollary~\ref{COR:smooth-locus-moment-calculations}.
Let $G,d,m,E_{\bm{c}}(q)$ be as in Problem~\ref{PROB:moments-of-failures-of-sqrt-cancellation}.
Corollary~\ref{COR:smooth-locus-moment-calculations} immediately follows from Corollary~\ref{COR:singular-locus-moment-calculations} (to Theorem~\ref{THM:main-general-result-informal}) and Proposition~\ref{PROP:full-locus-moment-calculations} (proven by elementary calculation) below.

\begin{corollary}
\label{COR:singular-locus-moment-calculations}
Fix $G$ with $d,m\geq 3$.
Uniformly over $p$, the following hold:
\begin{enumerate}
    \item $\EE_{\bm{c}\in \FF_p^m}[\abs{E_{\bm{c}}(p)}\cdot \bm{1}_{p\mid \disc(G,\bm{c})}]
    = O(p^{-1})\cdot p^{(m-3)/2}$.
    
    \item $\EE_{\bm{c}\in \FF_p^m}[\abs{E_{\bm{c}}(p^2)}\cdot \bm{1}_{p\mid \disc(G,\bm{c})}]
    = O(p^{-1})\cdot p^{m-3}$.
    
    \item $\EE_{\bm{c}\in \FF_p^m}[\abs{E_{\bm{c}}(p)}^2\cdot \bm{1}_{p\mid \disc(G,\bm{c})}]
    = O(p^{-1})\cdot p^{m-3}$.
\end{enumerate}
\end{corollary}

\begin{proposition}
\label{PROP:full-locus-moment-calculations}
Fix $G$ with $d\geq 2$ and $m\geq 3$.
Uniformly over $p$, the following hold:
\begin{enumerate}
    \item $\EE_{\bm{c}\in \FF_p^m}[E_{\bm{c}}(p)]
    = p^{-1}\cdot E(V_{\PP^{m-1}}(G)_{/\FF_p}) + O(p^{-1})$.
    
    \item $\EE_{\bm{c}\in \FF_p^m}[E_{\bm{c}}(p^2)]
    = (1+O(p^{-1/2}\bm{1}_{m=3})+O(p^{-1}))\cdot p^{m-3}$.
    
    \item $\EE_{\bm{c}\in \FF_p^m}[E_{\bm{c}}(p)^2]
    = (1+O(p^{-1/2}\bm{1}_{m=3})+O(p^{-1}))\cdot p^{m-3}$.
\end{enumerate}
\end{proposition}

For the proofs,
let $W\defeq V_{\PP^{m-1}}(G)_{/\FF_p}$ (given $p$) and $W_{\bm{c}}\defeq V_{\PP^{m-1}}(G,\bm{c}\cdot\bm{x})_{/\FF_p}$ (given $\bm{c}\in \FF_p^m$), and let $\disc(G)\in \ZZ$ denote the discriminant of $G$ (up to sign).

\begin{proof}
[Proof of Corollary~\ref{COR:singular-locus-moment-calculations}]
It suffices to prove the estimates for $p\nmid d(d-1)\disc(G)$.
Given such a prime $p$, we now make three observations:
\begin{itemize}
    \item By the definitions, $E_{\bm{0}}(q) = E(W_{\FF_q})+q^{m-2} \ll q^{m-2}$.
    So $p^{-m}E_{\bm{0}}(p)\ll p^{-2}$ and $p^{-m}E_{\bm{0}}(p)^2, p^{-m}E_{\bm{0}}(p^2)\ll p^{m-4}$.
    
    \item By Theorem~\ref{THM:Zak's-principle}, Proposition~\ref{PROP:amplification-for-projective-varieties}, and \cite{hooley1991number}*{Theorem~2},
    \begin{equation*}
    E_{\bm{c}}(q)/q^{(m-3)/2} \ll 1
    + q^{1/2}\cdot \bm{1}_{W_{\bm{c}}/\FF_p\;\textnormal{is $\abs{E}$-bad}}
    \end{equation*}
    holds for all $(\bm{c},q)\in(\FF_p^m\setminus \set{\bm{0}})\times \set{p, p^2, p^3, \ldots}$.
    
    \item By Theorem~\ref{THM:main-general-result-informal}, and Lang--Weil for $V(\disc(G,-))_{\map{sing}}$, we have
    \begin{equation*}
    \EE_{\bm{c}\in \FF_p^m}[\bm{1}_{p\mid \disc(G,\bm{c})}
    \cdot \bm{1}_{p\nmid \bm{c}}
    \cdot \bm{1}_{W_{\bm{c}}/\FF_p\;\textnormal{is $\abs{E}$-bad}}]
    \ll p^{-2}.
    \end{equation*}
\end{itemize}
By Lang--Weil for $V(\disc(G,-))$, the following bounds follow:
\begin{enumerate}
    \item $\EE_{\bm{c}\in \FF_p^m}[\abs{E_{\bm{c}}(p)}\cdot \bm{1}_{p\mid \disc(G,\bm{c})}]
    \ll p^{-2} + (1\cdot p^{-1} + p^{1/2}\cdot p^{-2})p^{(m-3)/2}$.
    
    \item $\EE_{\bm{c}\in \FF_p^m}[\abs{E_{\bm{c}}(p^2)}\cdot \bm{1}_{p\mid \disc(G,\bm{c})}]
    \ll p^{m-4} + (1\cdot p^{-1} + (p^2)^{1/2}\cdot p^{-2})p^{m-3}$.
    
    \item $\EE_{\bm{c}\in \FF_p^m}[\abs{E_{\bm{c}}(p)}^2\cdot \bm{1}_{p\mid \disc(G,\bm{c})}]
    \ll p^{m-4} + (1\cdot p^{-1} + (p^{1/2})^2\cdot p^{-2})p^{m-3}$.
\end{enumerate}
The result follows.
\end{proof}

In what follows, we write $\bm{c}\perp \bm{x}$ if $\bm{c}\cdot \bm{x} = 0$.
(Here $\bm{c}\cdot \bm{x}\defeq c_1x_1+\dots+c_mx_m$.)

\begin{proof}
[Proof of Proposition~\ref{PROP:full-locus-moment-calculations}]
Let $m_\ast\defeq m-3$.
We may restrict attention to $p\nmid \disc(G)$.
Let $p$ denote such a prime.
If $[\bm{x}]\in W(\FF_p)$, then $\PP_{\bm{c}\in \FF_p^m}[\bm{c}\perp \bm{x}] = p^{-1}$; so
\begin{equation}
\label{EQN:W_c-first-moment-identity}
\EE_{\bm{c}\in\FF_p^m}[\#W_{\bm{c}}(\FF_p)] = p^{-1}\cdot \card{W(\FF_p)}.
\end{equation}
By the definitions, $\#W_{\bm{c}}(\FF_p) = \card{\PP^{m_\ast}(\FF_p)} + E_{\bm{c}}(p)$ and $\card{W(\FF_p)} = \card{\PP^{1+m_\ast}(\FF_p)} + E(W)$.
But $p^{-1}\card{\PP^{1+m_\ast}(\FF_p)} - \card{\PP^{m_\ast}(\FF_p)} = p^{-1}$.
So Proposition~\ref{PROP:full-locus-moment-calculations}(1) follows from eq.~\eqref{EQN:W_c-first-moment-identity}.

Via $E_{\bm{c}}(p) = \card{W_{\bm{c}}(\FF_p)} - \card{\PP^{m_\ast}(\FF_p)}$, eq.~\eqref{EQN:W_c-first-moment-identity} also implies
\begin{equation}
\label{EQN:E_c-second-moment-opening}
\EE_{\bm{c}\in \FF_p^m}[E_{\bm{c}}(p)^2]
= \EE_{\bm{c}\in \FF_p^m}[\card{W_{\bm{c}}(\FF_p)}^2]
- 2(p^{-1}\card{W(\FF_p)})\cdot \card{\PP^{m_\ast}(\FF_p)}
+ \card{\PP^{m_\ast}(\FF_p)}^2.
\end{equation}
Meanwhile, by counting hyperplanes through pairs $([\bm{x}],[\bm{y}])\in W(\FF_p)^2$, we get
\begin{equation}
\label{EQN:W_c-second-moment-identity}
\EE_{\bm{c}\in \FF_p^m}[\card{W_{\bm{c}}(\FF_p)}^2]
= p^{-1}\card{W(\FF_p)}
+ p^{-2}(\card{W(\FF_p)}^2-\card{W(\FF_p)}).
\end{equation}
Upon plugging eq.~\eqref{EQN:W_c-second-moment-identity} into eq.~\eqref{EQN:E_c-second-moment-opening}, we obtain
\begin{equation*}
\EE_{\bm{c}\in \FF_p^m}[E_{\bm{c}}(p)^2]
= (p^{-1}-p^{-2})\card{W(\FF_p)}
+ (p^{-1}\card{W(\FF_p)} - \card{\PP^{m_\ast}(\FF_p)})^2.
\end{equation*}
But the Weil conjectures furnish the estimates
$(p^{-1}-p^{-2})\card{W(\FF_p)} = (p^{m_\ast}-p^{-2}) + O(p^{(m_\ast-1)/2})$ and $(p^{-1}\card{W(\FF_p)} - \card{\PP^{m_\ast}(\FF_p)})^2 \ll (p^{(m_\ast-1)/2})^2$.
So Proposition~\ref{PROP:full-locus-moment-calculations}(3) follows.

It remains to prove Proposition~\ref{PROP:full-locus-moment-calculations}(2).
Let $q\defeq p^2$.
If $[\bm{x}]\in W(\FF_q)$, then $\PP_{\bm{c}\in \FF_p^m}[\bm{c}\perp \bm{x}] = p^{-\deg \FF_p([\bm{x}])/\FF_p}$ (by Lemma~\ref{LEM:algebraic-vs-linear-dimension-coincidence} below).
It follows that
\begin{equation*}
\EE_{\bm{c}\in\FF_p^m}[\#W_{\bm{c}}(\FF_q)]
= p^{-1}\card{W(\FF_p)}
+ p^{-2}(\card{W(\FF_q)}-\card{W(\FF_p)}).
\end{equation*}
But by the Weil conjectures,
$(p^{-1}-p^{-2})\card{W(\FF_p)} = (p^{m_\ast}-p^{-2}) + O(p^{(m_\ast-1)/2})$ and $p^{-2}\card{W(\FF_q)} = \card{\PP^{m_\ast}(\FF_q)}+q^{-1}+O(q^{(m_\ast-1)/2})$.
This suffices.
\end{proof}

\begin{lemma}
\label{LEM:algebraic-vs-linear-dimension-coincidence}
Fix an integer $m\geq 1$.
Let $L/K$ be a quadratic extension of fields.
Fix $\bm{x}\in L^m\setminus \set{\bm{0}}$.
Let $f\maps K^m\to L$ denote the $K$-linear map $\bm{c}\mapsto \bm{x}\cdot \bm{c}$.
Then the following are equivalent:
(1) $\rank{f} = 2$;
(2) $\im{f} = L$;
and (3) $\deg K([\bm{x}])/K = 2$, i.e.~$[\bm{x}]\in \PP^{m-1}(\ol{K})$ lies over a closed point $x\in \PP^{m-1}_K$ of degree $2$.
(Concretely here, if say $\bm{x} = (x_1,\dots,x_m)$ with $x_m\neq 0$, then $K([\bm{x}])$ denotes the field extension of $K$ generated by the ratios $x_1/x_m, \dots, x_{m-1}/x_m$.)
\end{lemma}

\begin{proof}
Routine.
(But this result does not generalize in the obvious way!)
\end{proof}




\section*{Acknowledgements}

I am grateful to Nick Katz, Remke Kloosterman, and Will Sawin for helpful comments on earlier versions of this work.
I also thank Brendan Hassett for generous help on the topic of cubic scrolls in cubic fourfolds,
Peter Sarnak for helpful general comments,
and Will Sawin for sharing his knowledge of the closest existing literature and techniques.
I thank Anshul Adve and Liyang Yang for some uplifting preliminary discussions.
Finally, I thank the managing editors and the anonymous referee for their useful comments, especially on writing,
and also thank the anonymous referee for sharing some very nice geometric insights and remarks.

\bibliographystyle{amsxport}
\bibliography{final.bib}

\begin{bibdiv}
\begin{biblist}

\bib{benoist2012degres}{article}{
      author={Benoist, Olivier},
       title={Degr\'{e}s d'homog\'{e}n\'{e}it\'{e} de l'ensemble des
  intersections compl\`etes singuli\`eres},
        date={2012},
        ISSN={0373-0956},
     journal={Ann. Inst. Fourier (Grenoble)},
      volume={62},
      number={3},
       pages={1189\ndash 1214},
         url={https://doi.org/10.5802/aif.2720},
      review={\MR{3013820}},
}

\bib{MO211153intersect_components}{misc}{
      author={Bright, Martin},
       title={Intersection of two components and {J}acobian criterion},
organization={MathOverflow},
        date={2015},
        note={URL: \url{https://mathoverflow.net/q/211153} (version:
  2015-07-09)},
}

\bib{bombieri1967local}{article}{
      author={Bombieri, E.},
      author={Swinnerton-Dyer, Peter},
       title={On the local zeta function of a cubic threefold},
        date={1967},
        ISSN={0391-173X},
     journal={Ann. Scuola Norm. Sup. Pisa Cl. Sci. (3)},
      volume={21},
       pages={1\ndash 29},
      review={\MR{212019}},
}

\bib{carlini2008complete}{incollection}{
      author={Carlini, Enrico},
      author={Chiantini, Luca},
      author={Geramita, Anthony~V.},
       title={Complete intersections on general hypersurfaces},
        date={2008},
      volume={57},
       pages={121\ndash 136},
         url={https://doi.org/10.1307/mmj/1220879400},
        note={Special volume in honor of Melvin Hochster},
      review={\MR{2492444}},
}

\bib{DaoMO49299product_ideal}{misc}{
      author={Dao, Hailong},
       title={When is the product of two ideals equal to their intersection?},
organization={MathOverflow},
        date={2022},
        note={URL: \url{https://mathoverflow.net/q/49299} (version:
  2022-01-04)},
}

\bib{debarre2003lines}{article}{
      author={Debarre, Olivier},
       title={Lines on smooth hypersurfaces},
        date={2003-11},
     journal={Preprint},
      eprint={https://www.math.ens.fr/~debarre/Lines_hypersurfaces.pdf},
}

\bib{dimca1990betti}{article}{
      author={Dimca, Alexandru},
       title={Betti numbers of hypersurfaces and defects of linear systems},
        date={1990},
        ISSN={0012-7094},
     journal={Duke Math. J.},
      volume={60},
      number={1},
       pages={285\ndash 298},
         url={https://doi.org/10.1215/S0012-7094-90-06010-7},
      review={\MR{1047124}},
}

\bib{deligne1973groupes}{book}{
      author={Deligne, Pierre},
      author={Katz, Nicholas~M.},
       title={Groupes de monodromie en g\'{e}om\'{e}trie alg\'{e}brique. {II}},
      series={Lecture Notes in Mathematics, Vol. 340},
   publisher={Springer-Verlag, Berlin-New York},
        date={1973},
        note={S\'{e}minaire de G\'{e}om\'{e}trie Alg\'{e}brique du Bois-Marie
  1967--1969 (SGA 7 II), Dirig\'{e} par P. Deligne et N. Katz},
      review={\MR{0354657}},
}

\bib{dolgachev2016corrado}{incollection}{
      author={Dolgachev, Igor},
       title={Corrado {S}egre and nodal cubic threefolds},
        date={2016},
   booktitle={From classical to modern algebraic geometry},
      series={Trends Hist. Sci.},
   publisher={Birkh\"{a}user/Springer, Cham},
       pages={429\ndash 450},
         url={https://doi.org/10.1007/978-3-319-32994-9_11},
      review={\MR{3776662}},
}

\bib{eisenbud1987varieties}{incollection}{
      author={Eisenbud, David},
      author={Harris, Joe},
       title={On varieties of minimal degree (a centennial account)},
        date={1987},
   booktitle={Algebraic geometry, {B}owdoin, 1985 ({B}runswick, {M}aine,
  1985)},
      series={Proc. Sympos. Pure Math.},
      volume={46},
   publisher={Amer. Math. Soc., Providence, RI},
       pages={3\ndash 13},
         url={https://doi.org/10.1090/pspum/046.1/927946},
      review={\MR{927946}},
}

\bib{gelfand2008discriminants}{book}{
      author={Gelfand, I.~M.},
      author={Kapranov, M.~M.},
      author={Zelevinsky, A.~V.},
       title={Discriminants, resultants, and multidimensional determinants},
      series={Mathematics: Theory \& Applications},
   publisher={Birkh\"{a}user Boston, Inc., Boston, MA},
        date={1994},
        ISBN={0-8176-3660-9},
         url={https://doi.org/10.1007/978-0-8176-4771-1},
      review={\MR{1264417}},
}

\bib{ghorpade2008etale}{incollection}{
      author={Ghorpade, Sudhir~R.},
      author={Lachaud, Gilles},
       title={\'{E}tale cohomology, {L}efschetz theorems and number of points
  of singular varieties over finite fields},
        date={2002},
      volume={2},
       pages={589\ndash 631},
         url={https://doi.org/10.17323/1609-4514-2002-2-3-589-631},
        note={Dedicated to Yuri I. Manin on the occasion of his 65th birthday.
  Refer to the (corrected, revised and updated) version of August 15, 2008:
  \url{https://arxiv.org/abs/0808.2169v1}},
      review={\MR{1988974}},
}

\bib{grothendieck1972groupes}{book}{
      author={Grothendieck, Alexander},
       title={Groupes de monodromie en g\'{e}om\'{e}trie alg\'{e}brique. {I}},
      series={Lecture Notes in Mathematics, Vol. 288},
   publisher={Springer-Verlag, Berlin-New York},
        date={1972},
        note={S\'{e}minaire de G\'{e}om\'{e}trie Alg\'{e}brique du Bois-Marie
  1967--1969 (SGA 7 I), Dirig\'{e} par A. Grothendieck. Avec la collaboration
  de M. Raynaud et D. S. Rim},
      review={\MR{0354656}},
}

\bib{grimmelt2021representation}{article}{
      author={Grimmelt, Lasse},
      author={Sawin, Will},
       title={Representation of squares by nonsingular cubic forms},
        date={2021},
        ISSN={0021-2172},
     journal={Israel J. Math.},
      volume={242},
      number={2},
       pages={501\ndash 547},
         url={https://doi.org/10.1007/s11856-021-2116-2},
      review={\MR{4282090}},
}

\bib{hassett1996special}{thesis}{
      author={Hassett, Brendan},
       title={Special cubic hypersurfaces of dimension four},
        type={Ph.D. Thesis, Harvard University},
        date={1996},
  url={http://gateway.proquest.com/openurl?url_ver=Z39.88-2004&rft_val_fmt=info:ofi/fmt:kev:mtx:dissertation&res_dat=xri:pqdiss&rft_dat=xri:pqdiss:9631502},
      review={\MR{2694332}},
}

\bib{hooley1986HasseWeil}{article}{
      author={Hooley, Christopher},
       title={On {W}aring's problem},
        date={1986},
        ISSN={0001-5962},
     journal={Acta Math.},
      volume={157},
      number={1-2},
       pages={49\ndash 97},
         url={https://doi.org/10.1007/BF02392591},
      review={\MR{857679}},
}

\bib{hooley1991number}{article}{
      author={Hooley, Christopher},
       title={On the number of points on a complete intersection over a finite
  field},
        date={1991},
        ISSN={0022-314X},
     journal={J. Number Theory},
      volume={38},
      number={3},
       pages={338\ndash 358},
         url={https://doi.org/10.1016/0022-314X(91)90023-5},
        note={With an appendix by Nicholas M. Katz},
      review={\MR{1114483}},
}

\bib{hassett2010flops}{article}{
      author={Hassett, Brendan},
      author={Tschinkel, Yuri},
       title={Flops on holomorphic symplectic fourfolds and determinantal cubic
  hypersurfaces},
        date={2010},
        ISSN={1474-7480},
     journal={J. Inst. Math. Jussieu},
      volume={9},
      number={1},
       pages={125\ndash 153},
         url={https://doi.org/10.1017/S1474748009000140},
      review={\MR{2576800}},
}

\bib{katz2001sums}{incollection}{
      author={Katz, Nicholas~M.},
       title={Sums of {B}etti numbers in arbitrary characteristic},
        date={2001},
      volume={7},
       pages={29\ndash 44},
         url={https://doi.org/10.1006/ffta.2000.0303},
        note={Dedicated to Professor Chao Ko on the occasion of his 90th
  birthday},
      review={\MR{1803934}},
}

\bib{kloosterman2022maximal}{article}{
      author={Kloosterman, Remke},
       title={Maximal families of nodal varieties with defect},
        date={2022},
        ISSN={0025-5874},
     journal={Math. Z.},
      volume={300},
      number={2},
       pages={1141\ndash 1156},
         url={https://doi.org/10.1007/s00209-021-02814-7},
      review={\MR{4363772}},
}

\bib{kiehl2001weil}{book}{
      author={Kiehl, Reinhardt},
      author={Weissauer, Rainer},
       title={Weil conjectures, perverse sheaves and {$l$}'adic {F}ourier
  transform},
      series={Ergebnisse der Mathematik und ihrer Grenzgebiete. 3. Folge. A
  Series of Modern Surveys in Mathematics [Results in Mathematics and Related
  Areas. 3rd Series. A Series of Modern Surveys in Mathematics]},
   publisher={Springer-Verlag, Berlin},
        date={2001},
      volume={42},
        ISBN={3-540-41457-6},
         url={https://doi.org/10.1007/978-3-662-04576-3},
      review={\MR{1855066}},
}

\bib{lindner2020hypersurfaces}{article}{
      author={Lindner, Niels},
       title={Hypersurfaces with defect},
        date={2020},
        ISSN={0021-8693},
     journal={J. Algebra},
      volume={555},
       pages={1\ndash 35},
         url={https://doi.org/10.1016/j.jalgebra.2020.02.022},
      review={\MR{4081491}},
}

\bib{manin1968correspondences}{article}{
      author={Manin, Ju.~I.},
       title={Correspondences, motifs and monoidal transformations},
        date={1968},
     journal={Mat. Sb. (N.S.)},
      volume={77 (119)},
       pages={475\ndash 507},
      review={\MR{0258836}},
}

\bib{mustata2017lecture2}{misc}{
      author={Musta\c{t}\u{a}, Mircea},
       title={Lecture 2. {Cubic hypersurfaces I}: {Rationality} of certain
  cubic hypersurfaces},
organization={University of Michigan},
        date={2017},
        note={URL:
  \url{http://www-personal.umich.edu/~mmustata/lecture2_rationality.pdf}
  (accessed 2022-03-29)},
}

\bib{poonen2017rational}{book}{
      author={Poonen, Bjorn},
       title={Rational points on varieties},
      series={Graduate Studies in Mathematics},
   publisher={American Mathematical Society, Providence, RI},
        date={2017},
      volume={186},
        ISBN={978-1-4704-3773-2},
         url={https://doi.org/10.1090/gsm/186},
      review={\MR{3729254}},
}

\bib{poonen2020valuation}{article}{
      author={Poonen, Bjorn},
      author={Stoll, Michael},
       title={The valuation of the discriminant of a hypersurface},
        date={2020-08-26},
     journal={Preprint},
      eprint={https://math.mit.edu/~poonen/papers/discriminant.pdf},
}

\bib{saito2012discriminant}{article}{
      author={Saito, Takeshi},
       title={The discriminant and the determinant of a hypersurface of even
  dimension},
        date={2012},
        ISSN={1073-2780},
     journal={Math. Res. Lett.},
      volume={19},
      number={4},
       pages={855\ndash 871},
         url={https://doi.org/10.4310/MRL.2012.v19.n4.a10},
      review={\MR{3008420}},
}

\bib{serre2000local}{book}{
      author={Serre, Jean-Pierre},
       title={Local algebra},
      series={Springer Monographs in Mathematics},
   publisher={Springer-Verlag, Berlin},
        date={2000},
        ISBN={3-540-66641-9},
         url={https://doi.org/10.1007/978-3-662-04203-8},
        note={Translated from the French by CheeWhye Chin and revised by the
  author},
      review={\MR{1771925}},
}

\bib{skorobogatov1992exponential}{article}{
      author={Skorobogatov, Alexei~N.},
       title={Exponential sums, the geometry of hyperplane sections, and some
  {D}iophantine problems},
        date={1992},
        ISSN={0021-2172},
     journal={Israel J. Math.},
      volume={80},
      number={3},
       pages={359\ndash 379},
         url={https://doi.org/10.1007/BF02808077},
      review={\MR{1202578}},
}

\bib{SpeyerMO49261product_ideal}{misc}{
      author={Speyer, David~E.},
       title={When is the product of two ideals equal to their intersection?},
organization={MathOverflow},
        date={2010},
        note={URL: \url{https://mathoverflow.net/q/49261} (version:
  2010-12-13)},
}

\bib{stacks-project}{misc}{
      author={{Stacks project authors}, The},
       title={The {Stacks project}},
        date={2022},
        note={URL: \url{https://stacks.math.columbia.edu} (accessed
  2022-03-30)},
}

\bib{terakado2018determinant}{article}{
      author={Terakado, Yasuhiro},
       title={The determinant and the discriminant of a complete intersection
  of even dimension},
        date={2018},
        ISSN={1073-2780},
     journal={Math. Res. Lett.},
      volume={25},
      number={1},
       pages={255\ndash 280},
         url={https://doi.org/10.4310/mrl.2018.v25.n1.a11},
      review={\MR{3818622}},
}

\bib{github-singular-cubic-threefold-2021}{misc}{
      author={Wang, Victor~Y.},
       title={Tools for analyzing cubic threefolds over finite fields},
   publisher={GitHub},
        date={2021},
        note={URL:
  \url{https://github.com/wangyangvictor/singular_cubic_threefolds} (version:
  {\tt commit 5971de6})},
}

\bib{wang2022thesis}{thesis}{
      author={Wang, Victor~Y.},
       title={Families and dichotomies in the circle method},
        type={Ph.D. Thesis, Princeton University},
        date={2022},
        note={URL: \url{http://arks.princeton.edu/ark:/88435/dsp01rf55zb86g}},
}

\end{biblist}
\end{bibdiv}
\end{document}